\documentclass[final,3p,12pt]{elsarticle}
\usepackage{mathrsfs,enumerate,amssymb,amsmath,color,lineno}
\usepackage[all]{xy}
\usepackage{latexsym}
\usepackage{amsthm,a4wide,color}
\usepackage{amscd,verbatim}
\usepackage{graphicx}
\usepackage[colorlinks=true]{hyperref}
\usepackage{amssymb}
\theoremstyle{plain}
\newtheorem{theorem}{Theorem}

\theoremstyle{definition}
\newtheorem{definition}{Definition}
\newtheorem{example}{Example}

\numberwithin{equation}{section}

\newtheorem{remark}{Remark}

\journal{}

\newcommand{\pf}{\noindent \textbf{Proof.}\quad}
\newcommand{\epf}{\hspace{\stretch{1}}$\blacksquare$}
\begin{document}

\begin{frontmatter}
%% 标题
\title{Characterizing uninorms on bounded lattices\tnoteref{mytitlenote}}
%%基金资助信息
\tnotetext[mytitlenote]{This work was supported by the National Natural Science Foundation of China
(No. 11601449), the Science and Technology Innovation Team of Education Department of
Sichuan for Dynamical System and its Applications (No. 18TD0013), the Youth Science and Technology
Innovation Team of Southwest Petroleum University for Nonlinear Systems (No. 2017CXTD02), and the National
Natural Science Foundation of China (Key Program) (No. 51534006).}
%, and the National Nature Science Foundation of China (Key Program)
%(No. 51534006).}

%%第一作者
\author[a1,a2,a3]{Xinxing Wu\corref{mycorrespondingauthor}}
\cortext[mycorrespondingauthor]{Corresponding author}
\address[a1]{School of Sciences, Southwest Petroleum University, Chengdu, Sichuan 610500, China}
\address[a2]{Institute for Artificial Intelligence, Southwest Petroleum University, Chengdu, Sichuan 610500, China}
\address[a3]{Zhuhai College of Jilin University, Zhuhai, Guangdong 519041, China}
\ead{wuxinxing5201314@163.com}

\author[a4]{Guanrong Chen}
\address[a4]{Department of Electrical Engineering, City University of Hong Kong,
Hong Kong SAR, China}
\ead{gchen@ee.cityu.edu.hk}

%\author[a1]{Qin Zhang}
%\ead{2602558059@qq.com}
%
%\author[a4]{Xu Zhang\corref{mycorrespondingauthor}}
%\address[a4]{Department of Mathematics, Shandong University, Weihai, Shandong 264209, China}
%\ead{xu$\_$zhang$\_$sdu@mail.sdu.edu.cn}
%
%\author[a3]{Lidong Wang}
%\ead{wld0707@163.com}

%% 摘要
\begin{abstract}
This paper establishes  some equivalent conditions of a uninorm, extending an arbitrary triangular norm
on $[0, e]$ or an arbitrary triangular conorm on $[e, 1]$ to the whole lattice.
%%And, some illustrative examples are added for clarity.
\end{abstract}
%% 关键词
\begin{keyword}
Lattice; Uninorm; Incomparability.
\MSC[2010] 03B52, 06B20, 03E72.
\end{keyword}

\end{frontmatter}

\section{Introduction and preliminaries}%% 第一节 section
%Triangular norms ($t$-norms) and triangular conorms ($t$-conorms), being an extension of the conjunction and
%disjunction in classical two-valued logic, were systematically investigated by Schweizer and Sklar \cite{SS1983}.

As an extension of the logic connective conjunction and disjunction in classical two-valued logic,
triangular norms ($t$-norms) with the neutral $1$ and triangular conorms ($t$-conorms) with the neutral $0$
on the unit interval $[0, 1]$ were introduced by Menger \cite{Me1942} and by Schweizer and Sklar \cite{SS1961},
respectively. Then, the concepts were extended to uniorms by Yager and Rybalov~\cite{YR1996}, replacing the
neural element $1$ of the $t$-norm or neural element $0$ of the $t$-conorm by a value $e$ lying anywhere in $[0, 1]$.

On the other hand, a {\it lattice} \cite{Bir1967} is a partially ordered set $(L, \leq)$ satisfying that every pair of two
elements $x, y\in L$ have a greatest lower bound, called {\it infimum} and denoted as $x\wedge y$, as well as a smallest upper
bound, called {\it supremum} and denoted as $x\vee y$. A lattice is {\it bounded} if it has a top element and a
bottom element, written as $0$ and $1$, respectively. Let $(L,\leq,0,1)$ denote a bounded lattice with
top element $1$ and bottom element $0$ throughout this paper.

For $x, y\in L$, the expression $x<y$ means that $x\leq y$ and $x\neq y$. The elements $x$ and $y$ in $L$ are {\it comparable}
if $x\leq y$ or $y\leq x$. Otherwise, $x$ and $y$ are called {\it incomparable} if $x\nleq y$ and $y\nleq x$, and in this case,
it is denoted as $x\| y$. Let $I_e=\{x\in L: x \| e\}$ for $e\in L$.

%\begin{definition}{\rm \cite{Bir1967}
%Let $(L,~ \leq,~ 0,~ 1)$ be a bounded lattice. The elements $x$ and $y$ in $L$ are {\it comparable} if $x\leq y$ or
%$y\leq x$. Otherwise, $x$ and $y$ are called {\it incomparable} if $x\nleq y$ and $y\nleq x$.
%In this situation, we use $x\| y$ notation when $x$ and $y$ are incomparable.}
%\end{definition}
%
%Let

\begin{definition}\cite{Bir1967}
Let $(L,\leq,0,1)$ be a bounded lattice and $a, b\in L$ with $a\leq b$. Define a subinterval $[a, b]$ as
$$
[a, b]=\{x\in L: a\leq x\leq b\}.
$$
Other subintervals such as $[a, b)$ and $(a, b)$ are defined similarly.
Let $A(e)=\left([0, e]\times [e, 1]\right)\cup \left([e, 1]
\times [0, e]\right)$.
\end{definition}

\begin{definition}\cite{KM2015,YR1996}
Let $(L, \leq, 0, 1)$ be a bounded lattice. A binary operation $U: L^{2}\longrightarrow L$ is called
 a {\it uninorm on $L$} if it is commutative, associative, increasing with respect to both variables
 and there exist some element $e\in L$ called the {\it neutral element} such that $U(e, x)=x$ for all $x\in L$.
In particular, an operation $T: L^{2}\longrightarrow L$ is called a {\it $t$-norm} ({\it $t$-conorm})
if it is commutative, associative, increasing with respect
to both variables and has a neutral element $e=1$ ($e=0$).
\end{definition}

Uninorms on bounded lattices were originated from the work of Kara{\c{c}}al and Mesiar~\cite{KM2015}, whose
results guarantee the existence of uninorms with any neural element $e$ on any bounded
lattice $L$. Meanwhile, they also constructed the smallest uninorm and the greatest uninorm
with the neural element $e\in L$. Recently, constructing new uninorms on bounded lattices was investigated
by many authors (see \cite{BK2014,C2019,CK2017,CKM2016,DHQ2019,KEM2017,OZ2020,XL2020}).

%Let $(L, \leq, 0, 1)$ be a bounded lattice and $e\in L$. Let $A(e)= [0, e]\times [e, 1]\cup [e, 1]
%\times [0, e]$ and $I_e=\{x\in L: x \| e\}$. Let $A, B\subset L$.
%We say that $A$ and $B$ are {\it comparable} if there exist $x\in A$ and $y\in B$
%such that $x$ and $y$ are comparable. Otherwise, we say that $A$ and $B$ are {\it incomparable},
%i.e., $x\| y$ for all $x\in A$ and $y\in B$, in this case, we use the notation $A\| B$.

{\c{C}}ayl{\i} and Kara{\c{c}}al~\cite{CK2017} obtained the following results:

\begin{theorem}{\rm \cite[Theorem~3.9]{CK2017}}\label{CK-Thm-3.9}
Let $(X, \leq, 0, 1)$ be a bounded lattice and $e\in L\setminus \{0, 1\}$. Suppose that $x\vee y\in I_{e}$ for all
$x, y\in I_{e}$. If $T_{e}$ is a $t$-norm on $[0, e]$, then the function $U_{T}: L\times L\longrightarrow L$
defined by
\begin{equation}\label{U(T)-operation}
U_{T}(x, y)=
\begin{cases}
T_{e}(x, y), & (x, y)\in [0, e]^2, \\
1, & (x, y)\in (e, 1]^{2}, \\
y, & (x, y)\in [0, e]\times I_{e}, \\
x, & (x, y)\in I_{e}\times [0, e], \\
x \vee y, & \text{otherwise},
\end{cases}
\end{equation}
is a uninorm on $L$ with the neutral element $e$.
\end{theorem}

\begin{theorem}{\rm \cite[Theorem~3.12]{CK2017}}\label{CK-Thm-3.12}
Let $(X, \leq, 0, 1)$ be a bounded lattice and $e\in L\setminus \{0, 1\}$. Suppose that $x\wedge y\in I_{e}$ for all
$x, y\in I_{e}$. If $S_{e}$ is a $t$-conorm on $[e, 1]$, then the function $U_{S}: L\times L\longrightarrow L$
defined by
\begin{equation}\label{U(S)-operation}
U_{S}(x, y)=
\begin{cases}
0, & (x, y)\in [0, e)^2, \\
S_{e}(x, y), & (x, y)\in [e, 1]^{2}, \\
y, & (x, y)\in [0, e]\times I_{e}, \\
x, & (x, y)\in I_{e}\times [0, e], \\
x\wedge y, & \text{otherwise},
\end{cases}
\end{equation}
is a uninorm on $L$ with the neutral element $e$.
\end{theorem}

\begin{theorem}{\rm \cite[Theorem~3.1]{CK2017}}\label{CK-Thm-3.1}
Let $(X, \leq, 0, 1)$ be a bounded lattice and $e\in L\setminus \{0, 1\}$. Suppose that
either $x\vee y>e$ for all $x, y\in I_{e}$ or $x\vee y\in I_{e}$ for all $x, y\in I_{e}$.
If $T_{e}$ is a $t$-norm on $[0, e]$, then the function $U_{t}: L\times L\longrightarrow L$
defined by
\begin{equation}\label{U(t)-operation}
U_{t}(x, y)=
\begin{cases}
T_{e}(x, y), & (x, y)\in [0, e]^2, \\
x\vee y, & (x, y)\in A(e) \cup \left(I_{e}\times I_{e}\right), \\
y, & (x, y)\in [0, e]\times I_{e}, \\
x, & (x, y)\in I_{e}\times [0, e], \\
1, & \text{otherwise},
\end{cases}
\end{equation}
is a uninorm on $L$ with the neutral element $e$.
\end{theorem}

\begin{theorem}{\rm \cite[Theorem~3.5]{CK2017}}\label{CK-Thm-3.5}
Let $(X, \leq, 0, 1)$ be a bounded lattice and $e\in L\setminus \{0, 1\}$. Suppose that
either $x\wedge y<e$ for all $x, y\in I_{e}$ or $x\wedge y\in I_{e}$ for all $x, y\in I_{e}$.
If $S_{e}$ is a $t$-conorm on $[e, 1]$, then the function $U_{s}: L\times L\longrightarrow L$
defined by
\begin{equation}\label{U(s)-operation}
U_{s}(x, y)=
\begin{cases}
S_{e}(x, y), & (x, y)\in [e, 1]^2, \\
x\wedge y, & (x, y)\in A(e) \cup \left(I_{e}\times I_{e}\right), \\
y, & (x, y)\in [e, 1]\times I_{e}, \\
x, & (x, y)\in I_{e}\times [e, 1], \\
0, & \text{otherwise},
\end{cases}
\end{equation}
is a uninorm on $L$ with the neutral element $e$.
\end{theorem}

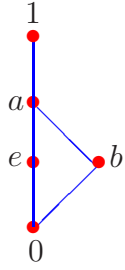
\begin{figure}[h]
\begin{center}
\begin{picture}(150,100)(100,-60)
\put(168,0){$a$}
\put(174.5,0){\color{red}$\bullet$}

\put(174.5,-22.4){\color{red}$\bullet$}
\put(168,-21){$e$}

\put(179,1){\color{blue}\line(1,-1){22}}
\put(199,-22.4){\color{red}$\bullet$}
\put(206,-21.5){$b$}

%\put(154,-21.5){\color{blue}\line(1,-1){24}}
\put(201,-21.5){\color{blue}\line(-1,-1){24}}
%\put(175,-28){\color{blue}\line(0,0){24}}
\put(174.5,-47){\color{red}$\bullet$}
\put(175.5,-57){$0$}

\put(177.5,-44){\color{blue}\line(0,0){72}}
\put(174.5,25){\color{red}$\bullet$}
\put(174.5,33){$1$}
\end{picture}
\end{center}
\caption{The Hasse diagram of the lattice $L_1$}
\label{Hasse-4}
\end{figure}

For $e\in L\setminus \{0, 1\}$, it is clear that $(0, e)\in [0, e]^{2} \cap A(e)$ and $(e, 1)\in [e, 1]^{2}
\cap A(e)$. According to the method of {\c{C}}ayl{\i} and Kara{\c{c}}al~\cite{CK2017} in proving the above Theorems~\ref{CK-Thm-3.1}
and \ref{CK-Thm-3.5}, it follows that
$$
0=T_{e}(0, e)=U_{t}(0, e)=0\vee e=e,
$$
and
$$
1=S_{e}(e, 1)=U_{s}(e, 1)=e\wedge 1=e,
$$
are both impossible. Meanwhile, consider the complete lattice $L_1$ with the Hasse diagram shown in Fig.~\ref{Hasse-4}
and take $S_{e}: \{e, a, 1\}^2 \longrightarrow \{e, a, 1\}$ by $S_{e}(x, y)=x\vee y$ for $(x, y)\in \{e, a, 1\}^{2}$.
Clearly, $I_{e}=\{b\}$. This implies that $x\wedge y=b\in I_{e}$ for all $x, y\in I_{e}$, i.e., $e$ satisfies the hypothesis
in Theorem~\ref{CK-Thm-3.12}. Applying Theorem~\ref{CK-Thm-3.12} (i.e., \cite[Theorem~3.12]{CK2017}) yields that $U_S$
defined by \eqref{U(S)-operation} is a uninorm on $L_1$. However, $U_{S}(0, a)=0<b=U_{S}(0, b)$, which is impossible, because
$a>b$. This shows that \cite[Theorem~3.12]{CK2017} does not hold.
%
%\begin{table}[h]
%\centering
%\begin{tabular}{|l|lllll|}
%\hline
% $U_s$& $0$ & $e$ & $a$ & $1$ & $b$ \\
% \hline
% ~$0$ & $0$ & $0$ & $b$ & $c$ & $0$ \\
% ~$e$ & $0$ & $a$ & $b$ & $c$ & $a$ \\
% ~$a$ & $b$ & $b$ & $b$ & $1$ & $b$ \\
% ~$1$ & $c$ & $c$ & $1$ & $c$ & $c$ \\
% ~$b$ & $0$ & $a$ & $b$ & $c$ & $e$ \\
% \hline
%\end{tabular}
%\label{tab-b}
%\caption{The binary operation $U_s$}
%\end{table}
Thus, we modify the above Theorems~\ref{CK-Thm-3.12}, \ref{CK-Thm-3.1},
and \ref{CK-Thm-3.5} in the following (see Theorems \ref{U(S)-operation-Thm}, \ref{U(t)-operation-Thm},
and \ref{U(s)-operation-Thm}) and obtain some equivalent characterizations of the binary operation defined
by \ref{U(T)-operation} to be uninorms. Similarly, we also give some equivalent characterizations of
the binary operations constructed by A{\c{s}}{\i}c{\i} and Mesiar~\cite{AM2020} to be uninorms.
These extend the main results in \cite{AM2020,CK2017}.

\section{Modification of Theorem~\ref{CK-Thm-3.12}}\label{S-2}
This section is devoted to correcting Theorem~\ref{CK-Thm-3.12} above and characterizing
uninorm property for the operation induced by a conorm on $[e, 1]$ (see Theorem~\ref{U(S)-operation-Thm} below).
In particular, our result shows that this operation is a uninorm with the neural element $e$ if and only if
$I_{e}\wedge I_{e}\subset I_{e}\cup \{0\}$, implying that the uninorm property of this operation is only
related to the property of $e$ in the lattice $L$.

\begin{theorem}\label{U(S)-operation-Thm}
Let $(X, \leq, 0, 1)$ be a bounded lattice and $e\in L\setminus \{0, 1\}$.
If $S_{e}$ is a $t$-conorm on $[e, 1]$ and the function $U_{S}: L\times L\longrightarrow L$
is defined by
\begin{equation}\label{U(S)-operation*}
U_{S}(x, y)=
\begin{cases}
0, & (x, y)\in [0, e)^2, \\
S_{e}(x, y), & (x, y)\in [e, 1]^2, \\
y, & (x, y)\in [e, 1]\times I_{e}, \\
x, & (x, y)\in I_{e}\times [e, 1], \\
x\wedge y, & \text{otherwise},
\end{cases}
\end{equation}
then the following statements are equivalent:
\begin{enumerate}[{\rm (i)}]
\item\label{4.i} $U_{S}$ is a uninorm on $L$ with the neutral element $e$;

\item\label{4.ii} $I_{e}=\emptyset$; otherwise, for any $y, z\in I_{e}$, it holds that $y\wedge z\in I_{e}$ or $y\wedge z=0$;

\item\label{4.iii} $I_{e}\wedge I_{e}\subset I_{e}\cup \{0\}$;

\item\label{4.iv} $\wedge$ is a $t$-norm on $I_{e} \cup \{0, 1\}$.
\end{enumerate}
%%$I_{e} \vee I_{e} \subset I_{e}$ or $I_{e}\vee I_{e}=\{1\}$.
\end{theorem}

\begin{figure}[h]
\begin{center}
\scalebox{0.6}{\includegraphics{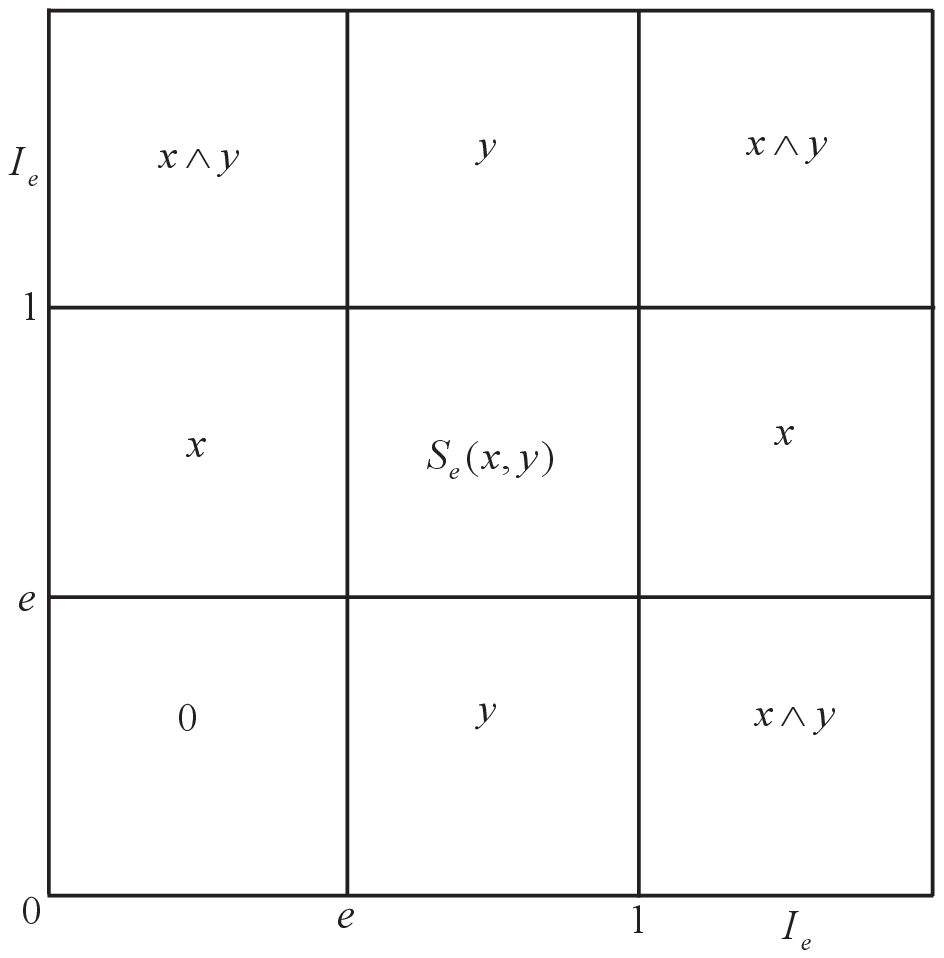}}
\renewcommand{\figure}{Fig.}
\caption{Structure of $U_{S}$ in Theorem~\ref{U(S)-operation-Thm}.
}
\end{center}
\end{figure}

\pf
Clearly, (\ref{4.ii})$\Longleftrightarrow$(\ref{4.iii})$\Longleftrightarrow$(\ref{4.iv}). It remains to show that
(\ref{4.i})$\Longleftrightarrow$(\ref{4.ii}).

\medskip

(\ref{4.i})$\Longrightarrow$(\ref{4.ii}).
Suppose, on the contrary, that there exist $y_0, z_0\in I_{e}$ such that $y_0\wedge z_0\notin I_{e}$
and $y_0\wedge z_0>0$. This implies that $y_0\wedge z_0\in (0, e)$.
Then,
$$
U_{S}(y_0\wedge z_0, U_{S}(y_0, z_0))=U_{S}(y_0\wedge z_0, y_0\wedge z_0)=0,
$$
and
$$
U_S(U_S(y_0\wedge z_0, y_0), z_0)=U_{S}(y_0\wedge z_0, z_0)=y_0\wedge z_0>0,
$$
which contradict with the associativity of $U_{S}$.

\medskip

(\ref{ii})$\Longrightarrow$(\ref{i}).
Clearly, $U_{S}$ is a uninorm when $I_{e}=\emptyset$. Without loss of generality, assume
that $I_{e}\neq \emptyset$. Firstly, it is not difficult to verify
that $U_S$ is a commutative binary operation with the neutral element $e$.

i.1) {\it Monotonicity}. For any $x, y\in L$ with $x\leq y$, it holds that $U_S(x, z)\leq U_S(y, z)$ for all $z\in L$.
Consider the following cases:
\begin{itemize}
\item[1.] $x\in [0, e)$.
\begin{itemize}
\item[1.1.] If $y\in [0, e)$, then
\begin{equation}\label{eq-4.2*}
U_{S}(x, z)=
\begin{cases}
0, & z\in [0, e), \\
x, & z=e, \\
x, & z\in (e, 1], \\
x\wedge z, & z\in I_{e},
\end{cases}
\end{equation}
and
$$
U_{S}(y, z)=
\begin{cases}
0, & z\in [0, e), \\
y, & z=e, \\
y, & z\in (e, 1], \\
y\wedge z, & z\in I_{e},
\end{cases}
$$
implying that $U_{S}(x, z)\leq U_{S}(y, z)$.

\item[1.2.] If $y\in (e, 1]$, then,
$$
U_{S}(y, z)=
\begin{cases}
z, & z\in [0, e), \\
y, & z=e, \\
S_{e}(y, z), & z\in (e, 1], \\
z, & z\in I_{e}.
\end{cases}
$$
This, together with \eqref{eq-4.2*} and $S_{e}(y, z)\geq y$, implies that $U_S(x, z)\leq U_S(y, z)$.

\item[1.3.] If $y\in I_{e}$, then
\begin{equation}\label{eq-4.5*}
U_{S}(y, z)=
\begin{cases}
y\wedge z, & z\in [0, e), \\
y, & z=e, \\
y, & z\in (e, 1], \\
y\wedge z, & z\in I_{e}.
\end{cases}
\end{equation}
This, together with \eqref{eq-4.2*}, implies that $U_S(x, z)\leq U_S(y, z)$.

\item[1.4.] If $y=e$, then $U_{S}(y, z)=z$. This, together with \eqref{eq-4.2*},
implies that $U_S(x, z)\leq U_S(y, z)$.
\end{itemize}

\item[2.] $x\in (e, 1]$. From $x\leq y$, it follows that $y\in (e, 1]$. Then,
\begin{equation}\label{eq-4.3*}
U_S(x, z)=
\begin{cases}
z, & z\in [0, e), \\
x, & z=e, \\
S_{e}(x, z), & z\in (e, 1], \\
z, & z\in I_{e},
\end{cases}
\end{equation}
and
\begin{equation}\label{eq-4.4*}
U_S(y, z)=
\begin{cases}
z, & z\in [0, e), \\
y, & z=e, \\
S_{e}(y, z), & z\in (e, 1], \\
z, & z\in I_{e},
\end{cases}
\end{equation}
implying that $U_S(x, z)\leq U_S(y, z)$, because $T_{e}$ is increasing.

\item[3.] $x=e$. This implies that $y\in [e, 1]$. Clearly, $U_S(x, z)=z= U_S(y, z)$ if $y=e$. If $y\in (e, 1]$,
applying \eqref{eq-4.4*} and $U_{S}(x, z)=z$, it follows that $U_{S}(x, z)\leq U_{S}(y, z)$.

\item[4.] $x\in I_{e}$. From $x\leq y$, it follows that $y\in (e, 1] \cup I_{e}$.
Applying \eqref{U(S)-operation*} yields that
\begin{equation}\label{4-10}
U_{S}(x, z)=
\begin{cases}
x\wedge z, & z\in [0, e), \\
x, & z=e, \\
x, & z\in (e, 1], \\
x\wedge z, & z\in I_{e}.
\end{cases}
\end{equation}
\begin{itemize}
\item[4.1.] If $y\in (e, 1]$, from \eqref{eq-4.4*} and \eqref{4-10},
it follows that $U_{S}(x, z)\leq U_{S}(y, z)$.

\item[4.2.] If $y\in I_{e}$, from \eqref{eq-4.5*} and \eqref{4-10}, it is clear that $U_{S}(x, z)\leq U_{S}(y, z)$.
\end{itemize}
\end{itemize}

i.2) {\it Associativity.} For any $x, y, z\in L$, it holds that $U_S(x, U_S(y, z))=U_S(U_S(x, y),z)$.%% for any $x, y, z\in L$.

If one of the elements $x$, $y$ and $z$ is equal to $e$, then the equality is always satisfied.
Otherwise, consider the following cases:
\begin{itemize}
\item[1.] $x\in [0, e)$.
\begin{itemize}
\item[1.1.] $y\in [0, e)$.
\begin{itemize}
\item[1.1.1.] If $z\in [0, e)$, then $U_{S}(x, U_{S}(y, z))=0=U_{S}(U_{S}(x, y), z)$.

\item[1.1.2.] If $z\in (e, 1]$, then $U_{S}(x, U_{S}(y, z))=U_{S}(x, y)=0$ and $U_{S}(U_{S}(x, y), z)
=U_{S}(0, z)=0$, implying that $U_{S}(x, U_{S}(y, z))=U_{S}(U_{S}(x, y),z)$.

\item[1.1.3.] If $z\in I_{e}$, then $U_{S}(x, U_{S}(y, z))=U_{S}(x, y\wedge z)=0$ and
$U_{S}(U_{S}(x, y), z)=U_{S}(0, z)=0$, implying that $U_{S}(x, U_{S}(y, z))=U_{S}(U_{S}(x, y),z)$.
\end{itemize}
\item[1.2.] $y\in (e, 1]$.
\begin{itemize}
\item[1.2.1.] If $z\in [0, e)$, then $U_{S}(x, U_{S}(y, z))=U_{S}(x, z)=0$ and
$U_{S}(U_{S}(x, y), z)=U_{S}(x, z)=0$, implying that $U_{S}(x, U_{S}(y, z))=U_{S}(U_{S}(x, y),z)$.

\item[1.2.2.] If $z\in (e, 1]$, then $U_{S}(x, U_{S}(y, z))=U_{S}(x, S_{e}(y, z))=x$ and
$U_{S}(U_{S}(x, y), z)=U_{S}(x, z)=x$, implying that $U_{S}(x, U_{S}(y, z))=U_{S}(U_{S}(x, y),z)$.

\item[1.2.3.] If $z\in I_{e}$, then $U_{S}(x, U_{S}(y, z))=U_{S}(x, z)=x\wedge z$ and
$U_{S}(U_{S}(x, y), z)=U_{S}(x, z)=x\wedge z$, implying that $U_{S}(x, U_{S}(y, z))=U_{S}(U_{S}(x, y),z)$.
\end{itemize}
\item[1.3.] $y\in I_{e}$.
\begin{itemize}
\item[1.3.1.] If $z\in [0, e)$, then $U_{S}(x, U_{S}(y, z))=U_{S}(x, y\wedge z)=0$ and
$U_{S}(U_{S}(x, y), z)=U_{S}(x\wedge y, z)=0$, implying that $U_{S}(x, U_{S}(y, z))=U_{S}(U_{S}(x, y),z)$.

\item[1.3.2.] If $z\in (e, 1]$, then $U_{S}(x, U_{S}(y, z))=U_{S}(x, y)=x\wedge y$ and
$U_{S}(U_{S}(x, y), z)=U_{S}(x\wedge y, z)=x\wedge y$, implying that $U_{S}(x, U_{S}(y, z))=U_{S}(U_{S}(x, y),z)$.

\item[1.3.3.] If $z\in I_{e}$, then
\begin{align*}
U_{S}(x, U_{S}(y, z))
&=U_{S}(x, y\wedge z)\\
&=
\begin{cases}
x\wedge y\wedge z, & y\vee z\in I_{e}, \\
0, & y\wedge z=0,
\end{cases}\\
&=x\wedge y\wedge z,
\end{align*}
and
$U_{S}(U_{S}(x, y), z)=U_{S}(x\wedge y, z)=x\wedge y\wedge z$, implying that $U_{S}(x, U_{S}(y, z))=U_{S}(U_{S}(x, y),z)$.
\end{itemize}
\end{itemize}

\item[2.] $x\in (e, 1]$.
\begin{itemize}
\item[2.1.] $y\in [0, e)$.
\begin{itemize}
\item[2.1.1.] If $z\in [0, e)$, then
\begin{align*}
U_{S}(x, U_{S}(y, z))&=U_{S}(U_{S}(z, y), x) ~~ (\text{commutativity of } U_{S})\\
&=U_{S}(z, U_{S}(y, x))~~ (\text{by 1.1.2})\\
&=U_{S}(U_{S}(x, y), z)~~ (\text{commutativity of } U_{S}).
\end{align*}

\item[2.1.2.] If $z\in (e, 1]$, then $U_{S}(x, U_{S}(y, z))=U_{S}(x, z)$ and
$U_{S}(U_{S}(x, y), z)=U_{S}(x, z)$, implying that $U_{S}(x, U_{S}(y, z))=U_{S}(U_{S}(x, y),z)$.

\item[2.1.3.] If $z\in I_{e}$, then $U_{S}(x, U_{S}(y, z))=U_{S}(x, z)$ and
$U_{S}(U_{S}(x, y), z)=U_{S}(x, z)$, implying that $U_{S}(x, U_{S}(y, z))=U_{S}(U_{S}(x, y),z)$.
\end{itemize}
\item[2.2.] $y\in (e, 1]$.
\begin{itemize}
\item[2.2.1.] If $z\in [0, a_1)$, then
\begin{align*}
U_{S}(x, U_{S}(y, z))&= U_{S}(U_{S}(z, y), x)~~ (\text{commutativity of } U_{S})\\
&= U_{S}(z, U_{S}(y, x))~~ (\text{by 1.2.2})\\
&= U_{S}(U_{S}(x, y), z)~~ (\text{commutativity of } U_{S}).
\end{align*}

\item[2.2.2.] If $z\in (e, 1]$, then $U_{S}(x, U_{S}(y, z))=U_{S}(x, S_{e}(y, z))=S_{e}(x, S_{e}(y, z))
=S_{e}(S_{e}(x, y), z)=S_{e}(U_{S}(x, y), z)=U_{S}(U_{S}(x, y), z)$.

\item[2.2.3.] If $z\in I_{e}$, then $U_{S}(x, U_{S}(y, z))=U_{S}(x, z)=z$ and
$U_{S}(U_{S}(x, y), z)=U_{S}(S_{e}(x, y), z)=z$, implying that $U_{S}(x, U_{S}(y, z))=U_{S}(U_{S}(x, y),z)$.
\end{itemize}

\item[2.3.] $y\in I_{e}$.
\begin{itemize}
\item[2.3.1.] If $z\in [0, e)$, then
\begin{align*}
U_{S}(x, U_{S}(y, z))&=U_{S}(U_{S}(z, y), x) ~~ (\text{commutativity of } U_{S})\\
&=U_{S}(z, U_{S}(y, x))~~ (\text{by 1.3.2})\\
&=U_{S}(U_{S}(x, y), z)~~ (\text{commutativity of } U_{S}).
\end{align*}

\item[2.3.2.] If $z\in (e, 1]$, then $U_{S}(x, U_{S}(y, z))=U_{S}(x, y)=y$ and
$U_{S}(U_{S}(x, y), z)=U_{S}(y, z)=y$, implying that $U_{S}(x, U_{S}(y, z))=U_{S}(U_{S}(x, y),z)$.

\item[2.3.3.] If $z\in I_{e}$, then
\begin{align*}
U_{S}(x, U_{S}(y, z))&=U_{S}(x, y\vee z)\\
&=\begin{cases}
y\wedge z, & y\wedge z\in I_{e}, \\
0, & y\wedge z=0,
\end{cases}\\
& =y\wedge z,
\end{align*}
and
$U_{S}(U_{S}(x, y), z)=U_{S}(y, z)=y\wedge z$,
implying that $U_{S}(x, U_{S}(y, z))=U_{S}(U_{S}(x, y),z)$.
\end{itemize}

\item[3.] $x\in I_{e}$.
\begin{itemize}
\item[3.1.] $y\in [0, e)$.
\begin{itemize}
\item[3.1.1.] If $z\in [0, e)$, then
\begin{align*}
U_{S}(x, U_{S}(y, z))&=U_{S}(U_{S}(z, y), x) ~~ (\text{commutativity of } U_{S})\\
&=U_{S}(z, U_{S}(y, x))~~ (\text{by 1.1.3})\\
&=U_{S}(U_{S}(x, y), z)~~ (\text{commutativity of } U_{S}).
\end{align*}

\item[3.1.2.] If $z\in (e, 1]$, then
\begin{align*}
U_{S}(x, U_{S}(y, z))&=U_{S}(U_{S}(z, y), x) ~~ (\text{commutativity of } U_{S})\\
&=U_{S}(z, U_{S}(y, x))~~ (\text{by 2.1.3})\\
&=U_{S}(U_{S}(x, y), z)~~ (\text{commutativity of } U_{S}).
\end{align*}

\item[3.1.3.] If $z\in I_{e}$, then $U_{S}(x, U_{S}(y, z))=U_{S}(x, y\wedge z)=x\wedge y\wedge z$ and
$U_{S}(U_{S}(x, y), z)=U_{S}(x\wedge y, z)=x\wedge y\wedge z$, implying that $U_{S}(x, U_{S}(y, z))=U_{S}(U_{S}(x, y),z)$.
\end{itemize}
\item[3.2.] $y\in (e, 1]$.
\begin{itemize}
\item[3.2.1.] If $z\in [0, e)$, then
\begin{align*}
U_{S}(x, U_{S}(y, z))&=U_{S}(U_{S}(z, y), x)~~ (\text{commutativity of } U_{S})\\
&= U_{S}(z, U_{S}(y, x))~~ (\text{by 1.2.3})\\
&= U_{S}(U_{S}(x, y), z)~~ (\text{commutativity of } U_{S}).
\end{align*}

\item[3.2.2.] If $z\in (e, 1]$, then
\begin{align*}
U_{S}(x, U_{S}(y, z))&= U_{S}(U_{S}(z, y), x)~~ (\text{commutativity of } U_{S})\\
&= U_{S}(z, U_{S}(y, x))~~ (\text{by 2.2.3})\\
&= U_{S}(U_{S}(x, y), z)~~ (\text{commutativity of } U_{S}).
\end{align*}

\item[3.2.3.] If $z\in I_{e}$, then $U_{S}(x, U_{S}(y, z))=U_{S}(x, z)=x\wedge z$ and
$U_{S}(U_{S}(x, y), z)=U_{S}(x, z)=x\wedge z$, implying that $U_{S}(x, U_{S}(y, z))=U_{S}(U_{S}(x, y),z)$.
\end{itemize}

\item[3.3.] $y\in I_{e}$.
\begin{itemize}
\item[3.3.1.] If $z\in [0, e)$, then
\begin{align*}
U_{S}(x, U_{S}(y, z))&=U_{S}(U_{S}(z, y), x)~~ (\text{commutativity of } U_{S})\\
&= U_{S}(z, U_{S}(y, x))~~ (\text{by 1.3.3})\\
&= U_{S}(U_{S}(x, y), z)~~ (\text{commutativity of } U_{S}).
\end{align*}

\item[3.3.2.] If $z\in (e, 1]$, then
\begin{align*}
U_{S}(x, U_{S}(y, z))&= U_{S}(U_{S}(z, y), x)~~ (\text{commutativity of } U_{S})\\
&= U_{S}(z, U_{S}(y, x))~~ (\text{by 2.3.3})\\
&= U_{S}(U_{S}(x, y), z)~~ (\text{commutativity of } U_{S}).
\end{align*}

\item[3.3.3.] If $z\in I_{e}$, then
\begin{align*}
U_{S}(x, U_{S}(y, z))&=U_{S}(x, y\wedge z)\\
&=\begin{cases}
x\wedge y\wedge z, & y\wedge z\in I_{e}, \\
0, & y\wedge z=0,
\end{cases}\\
&=x\wedge y \wedge z,
\end{align*} and
\begin{align*}
U_{S}(U_{S}(x, y), z)&=U_{S}(x\wedge y, z)\\
&=\begin{cases}
x\wedge y\wedge z, & x\wedge y\in I_{e}, \\
0, & x\wedge y=0,
\end{cases}\\
&=x\wedge y \wedge z,
\end{align*}
implying that $U_{S}(x, U_{S}(y, z))=U_{S}(U_{S}(x, y),z)$.
\end{itemize}
\end{itemize}
\end{itemize}
\end{itemize}
\epf

\section{Modifications of Theorems~\ref{CK-Thm-3.1} and \ref{CK-Thm-3.5}}
Similarly to Section~\ref{S-2}, this section is devoted to correcting Theorems~\ref{CK-Thm-3.1} and
\ref{CK-Thm-3.5} and characterizing the uninorm property (see Theorems~\ref{U(S)-operation-Thm} and \ref{U(s)-operation-Thm})
for some new binary operations on the lattices.
\begin{theorem}\label{U(t)-operation-Thm}
Let $(X, \leq, 0, 1)$ be a bounded lattice and $e\in L\setminus \{0, 1\}$.
If $T_{e}$ is a $t$-norm on $[0, e]$ and the function $U_{t}: L\times L\longrightarrow L$
is defined by
\begin{equation}\label{U(t)-operation*}
U_{t}(x, y)=
\begin{cases}
T_{e}(x, y), & (x, y)\in [0, e]^2, \\
x\vee y, & (x, y)\in \left([0, e]\times (e, 1]\right) \cup \left((e, 1]\times [0, e]\right) \cup \left(I_{e}\times I_{e}\right), \\
y, & (x, y)\in [0, e]\times I_{e}, \\
x, & (x, y)\in I_{e}\times [0, e], \\
1, & \text{otherwise},
\end{cases}
\end{equation}
then the following statements are equivalent:
\begin{enumerate}[{\rm (i)}]
\item\label{i} $U_{t}$ is a uninorm on $L$ with the neutral element $e$;

\item\label{ii} $I_{e}=\emptyset$; otherwise, for any $y, z\in I_{e}$, it holds that $y\vee z\in I_{e}$ or $y\vee z=1$;

\item\label{iii} $I_{e}\vee I_{e}\subset I_{e}\cup \{1\}$;

\item\label{iv} $\vee$ is a $t$-conorm on $I_{e} \cup \{0, 1\}$.
\end{enumerate}
%%$I_{e} \vee I_{e} \subset I_{e}$ or $I_{e}\vee I_{e}=\{1\}$.
\end{theorem}

\begin{figure}[h]
\begin{center}
\scalebox{0.6}{\includegraphics{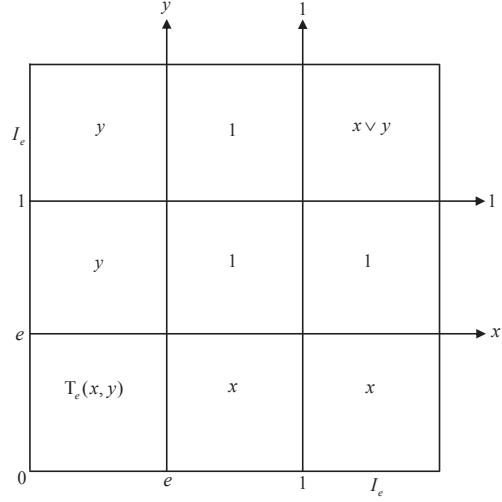}}
\renewcommand{\figure}{Fig.}
\caption{Structure of $U_{t}$ in Theorem~\ref{U(t)-operation-Thm}.
}
\end{center}
\end{figure}

\pf
Clearly, (\ref{ii})$\Longleftrightarrow$(\ref{iii})$\Longleftrightarrow$(\ref{iv}). It remains to show that
(\ref{i})$\Longleftrightarrow$(\ref{ii}).

\medskip

(\ref{i})$\Longrightarrow$(\ref{ii}).
Suppose, on the contrary, that there exist $y_0, z_0\in I_{e}$ such that $y_0\vee z_0\notin I_{e}$
and $y_0\vee z_0<1$. This implies that $y_0\vee z_0\in (e, 1)$.
Consequently,
$$
U_{t}(y_0, U_{t}(y_0, z_0))=U_{t}(y_0, y_0\vee z_0)=1,
$$
and
$$
U_t(U_t(y_0, y_0), z_0)=U_{t}(y_0, z_0)=y_0\vee z_0<1,
$$
which contradict with the associativity of $U_{t}$.

\medskip

(\ref{ii})$\Longrightarrow$(\ref{i}).
Clearly, $U_{t}$ is a uninorm when $I_{e}=\emptyset$. Without loss of generality, assume
that $I_{e}\neq \emptyset$. Firstly, it is not difficult to verify
that $U_t$ is a commutative binary operation with the neutral element $e$.

i.1) {\it Monotonicity}. For any $x, y\in L$ with $x\leq y$, it holds that $U_t(x, z)\leq U_t(y, z)$ for all $z\in L$.
Consider the following cases:
\begin{itemize}
\item[1.] $x\in [0, e)$.
\begin{itemize}
\item[1.1.] If $y\in [0, e)$, then
\begin{equation}\label{eq-2*}
U_{t}(x, z)=
\begin{cases}
T_{e}(x, z), & z\in [0, e), \\
x, & z=e, \\
z, & z\in (e, 1], \\
z, & z\in I_{e},
\end{cases}
\end{equation}
and
$$
U_{t}(y, z)=
\begin{cases}
T_{e}(y, z), & z\in [0, e), \\
y, & z=e, \\
z, & z\in (e, 1], \\
z, & z\in I_{e},
\end{cases}
$$
implying that $U_{t}(x, z)\leq U_{t}(y, z)$, because $T_{e}$ is increasing.

\item[1.2.] If $y\in (e, 1]$, then,
$$
U_{t}(y, z)=
\begin{cases}
y, & z\in [0, e), \\
y, & z=e, \\
1, & z\in (e, 1], \\
1, & z\in I_{e}.
\end{cases}
$$
This, together with \eqref{eq-2*} and $T_{e}(x, z)\leq x$, implies that $U_t(x, z)\leq U_t(y, z)$.

\item[1.3.] If $y\in I_{e}$, then
\begin{equation}\label{eq-5*}
U_{t}(y, z)=
\begin{cases}
y, & z\in [0, e), \\
y, & z=e, \\
1, & z\in (e, 1], \\
y\vee z, & z\in I_{e}.
\end{cases}
\end{equation}
This, together with \eqref{eq-2*} and $T_{e}(x, z)\leq x$, implies that $U_t(x, z)\leq U_t(y, z)$.

\item[1.4.] If $y=e$, then $U_{t}(y, z)=z$. This, together with \eqref{eq-2*}, implies that $U_t(x, z)\leq U_t(y, z)$.
\end{itemize}

\item[2.] $x\in (e, 1]$. From $x\leq y$, it follows that $y\in (e, 1]$. Then,
\begin{equation}\label{eq-3*}
U_t(x, z)=
\begin{cases}
x, & z\in [0, e), \\
x, & z=e, \\
1, & z\in (e, 1], \\
1, & z\in I_{e},
\end{cases}
\end{equation}
and
\begin{equation}\label{eq-4*}
U_t(y, z)=
\begin{cases}
y, & z\in [0, e), \\
y, & z=e, \\
1, & z\in (e, 1], \\
1, & z\in I_{e},
\end{cases}
\end{equation}
implying that $U_t(x, z)\leq U_t(y, z)$.

\item[3.] $x=e$. This implies that $y\in [e, 1]$. Clearly, $U_t(x, z)=z= U_t(y, z)$ if $y=e$. If $y\in (e, 1]$,
applying \eqref{eq-4*} and $U_{t}(x, z)=z$, it follows that $U_{t}(x, z)\leq U_{t}(y, z)$.

\item[4.] $x\in I_{e}$. From $x\leq y$, it follows that $y\in (e, 1] \cup I_{e}$.
Applying \eqref{U(t)-operation*} yields that
\begin{equation}\label{1-10}
U_{t}(x, z)=
\begin{cases}
x, & z\in [0, e), \\
x, & z=e, \\
1, & z\in (e, 1], \\
x\vee z, & z\in I_{e}.
\end{cases}
\end{equation}
\begin{itemize}
\item[4.1.] If $y\in (e, 1]$, from \eqref{eq-4*} and \eqref{1-10},
it follows that $U_{t}(x, z)\leq U_{t}(y, z)$.

\item[4.2.] If $y\in I_{e}$, from \eqref{eq-5*} and \eqref{1-10}, it is clear that $U_{t}(x, z)\leq U_{t}(y, z)$.
\end{itemize}
\end{itemize}

i.2) {\it Associativity.} For any $x, y, z\in L$, it holds that $U_t(x, U_t(y, z))=U_t(U_t(x, y),z)$.%% for any $x, y, z\in L$.

If one of the elements $x$, $y$ and $z$ is equal to $e$, then the equality is always satisfied.
Otherwise, consider the following cases:
\begin{itemize}
\item[1.] $x\in [0, e)$.
\begin{itemize}
\item[1.1.] $y\in [0, e)$.
\begin{itemize}
\item[1.1.1.] If $z\in [0, e)$, then $U_{t}(x, U_{t}(y, z))=U_{t}(x, T_e(y, z))=T_{e}(x, T_{e}(y, z))=
T_{e}(T_{e}(x, y), z)=T_{e}(U_{t}(x, y), z)=U_{t}(U_{t}(x, y), z)$.

\item[1.1.2.] If $z\in (e, 1]$, then $U_{t}(x, U_{t}(y, z))=U_{t}(x, z)=z$ and $U_{t}(U_{t}(x, y), z)
=U_{t}(T_{e}(x, y), z)=z$, implying that $U_t(x, U_t(y, z))=U_t(U_t(x, y),z)$.

\item[1.1.3.] If $z\in I_{e}$, then $U_{t}(x, U_{t}(y, z))=U_{t}(x, z)=z$ and
$U_{t}(U_{t}(x, y), z)=U_{t}(T_{e}(x, y), z)=z$, implying that $U_t(x, U_t(y, z))=U_t(U_t(x, y),z)$.
\end{itemize}
\item[1.2.] $y\in (e, 1]$.
\begin{itemize}
\item[1.2.1.] If $z\in [0, e)$, then $U_{t}(x, U_{t}(y, z))=U_{t}(x, y)=y$ and
$U_{t}(U_{t}(x, y), z)=U_{t}(y, z)=y$, implying that $U_t(x, U_t(y, z))=U_t(U_t(x, y),z)$.

\item[1.2.2.] If $z\in (e, 1]$, then $U_{t}(x, U_{t}(y, z))=U_{t}(x, 1)=1$ and
$U_{t}(U_{t}(x, y), z)=U_{t}(y, z)=1$, implying that $U_t(x, U_t(y, z))=U_t(U_t(x, y),z)$.

\item[1.2.3.] If $z\in I_{e}$, then $U_{t}(x, U_{t}(y, z))=U_{t}(x, 1)=1$ and
$U_{t}(U_{t}(x, y), z)=U_{t}(y, z)=1$, implying that $U_t(x, U_t(y, z))=U_t(U_t(x, y),z)$.
\end{itemize}
\item[1.3.] $y\in I_{e}$.
\begin{itemize}
\item[1.3.1.] If $z\in [0, e)$, then $U_{t}(x, U_{t}(y, z))=U_{t}(x, y)=y$ and
$U_{t}(U_{t}(x, y), z)=U_{t}(y, z)=y$, implying that $U_t(x, U_t(y, z))=U_t(U_t(x, y),z)$.

\item[1.3.2.] If $z\in (e, 1]$, then $U_{t}(x, U_{t}(y, z))=U_{T}(x, 1)=1$ and
$U_{t}(U_{t}(x, y), z)=U_{t}(y, z)=1$, implying that $U_t(x, U_t(y, z))=U_t(U_t(x, y),z)$.

\item[1.3.3.] If $z\in I_{e}$, then
\begin{align*}
U_{t}(x, U_{t}(y, z))
&=U_{t}(x, y\vee z)\\
&=
\begin{cases}
y\vee z, & y\vee z\in I_{e}, \\
1, & y\vee z=1,
\end{cases}\\
&=y\vee z,
\end{align*}
and
$U_{t}(U_{t}(x, y), z)=U_{t}(y, z)=y\vee z$, implying that $U_t(x, U_t(y, z))=U_t(U_t(x, y),z)$.
\end{itemize}
\end{itemize}

\item[2.] $x\in (e, 1]$.
\begin{itemize}
\item[2.1.] $y\in [0, e)$.
\begin{itemize}
\item[2.1.1.] If $z\in [0, e)$, then
\begin{align*}
U_t(x, U_t(y, z))&=U_t(U_t(z, y), x) ~~ (\text{commutativity of } U_t)\\
&=U_t(z, U_t(y, x))~~ (\text{by 1.1.2})\\
&=U_t(U_t(x, y), z)~~ (\text{commutativity of } U_t).
\end{align*}

\item[2.1.2.] If $z\in (e, 1]$, then $U_{t}(x, U_{t}(y, z))=U_{t}(x, z)$ and
$U_{t}(U_{t}(x, y), z)=U_{t}(x, z)$, implying that $U_t(x, U_t(y, z))=U_t(U_t(x, y),z)$.

\item[2.1.3.] If $z\in I_{e}$, then $U_{t}(x, U_{t}(y, z))=U_{t}(x, z)$ and
$U_{t}(U_{t}(x, y), z)=U_{t}(x, z)$, implying that $U_t(x, U_t(y, z))=U_t(U_t(x, y),z)$.
\end{itemize}
\item[2.2.] $y\in (e, 1]$.
\begin{itemize}
\item[2.2.1.] If $z\in [0, a_1)$, then
\begin{align*}
U_t(x, U_t(y, z))&= U_t(U_t(z, y), x)~~ (\text{commutativity of } U_t)\\
&= U_t(z, U_t(y, x))~~ (\text{by 1.2.2})\\
&= U_t(U_t(x, y), z)~~ (\text{commutativity of } U_t).
\end{align*}

\item[2.2.2.] If $z\in (e, 1]$, then $U_{t}(x, U_{t}(y, z))=U_{t}(x, 1)=1$ and
$U_{t}(U_{t}(x, y), z)=U_{t}(1, z)=1$, implying that $U_t(x, U_t(y, z))=U_t(U_t(x, y),z)$.

\item[2.2.3.] If $z\in I_{e}$, then $U_{t}(x, U_{t}(y, z))=U_{t}(x, 1)=1$ and
$U_{t}(U_{t}(x, y), z)=U_{t}(1, z)=1$, implying that $U_t(x, U_t(y, z))=U_t(U_t(x, y),z)$.
\end{itemize}

\item[2.3.] $y\in I_{e}$.
\begin{itemize}
\item[2.3.1.] If $z\in [0, e)$, then
\begin{align*}
U_t(x, U_t(y, z))&=U_t(U_t(z, y), x) ~~ (\text{commutativity of } U_t)\\
&=U_t(z, U_t(y, x))~~ (\text{by 1.3.2})\\
&=U_t(U_t(x, y), z)~~ (\text{commutativity of } U_t).
\end{align*}

\item[2.3.2.] If $z\in (e, 1]$, then $U_{t}(x, U_{t}(y, z))=U_{t}(x, 1)=1$ and
$U_{t}(U_{t}(x, y), z)=U_{t}(1, z)=1$, implying that $U_t(x, U_t(y, z))=U_t(U_t(x, y),z)$.

\item[2.3.3.] If $z\in I_{e}$, then
\begin{align*}
U_{t}(x, U_{t}(y, z))&=U_{t}(x, y\vee z)\\
&=\begin{cases}
1, & y\vee z\in I_{e}, \\
1, & y\vee z=1,
\end{cases}\\
& =1,
\end{align*}
and
$U_{t}(U_{t}(x, y), z)=U_{t}(1, z)=1$,
implying that $U_t(x, U_t(y, z))=U_t(U_t(x, y),z)$.
\end{itemize}

\item[3.] $x\in I_{e}$.
\begin{itemize}
\item[3.1.] $y\in [0, e)$.
\begin{itemize}
\item[3.1.1.] If $z\in [0, e)$, then
\begin{align*}
U_t(x, U_t(y, z))&=U_t(U_t(z, y), x) ~~ (\text{commutativity of } U_t)\\
&=U_t(z, U_t(y, x))~~ (\text{by 1.1.3})\\
&=U_t(U_t(x, y), z)~~ (\text{commutativity of } U_t).
\end{align*}

\item[3.1.2.] If $z\in (e, 1]$, then
\begin{align*}
U_t(x, U_t(y, z))&=U_t(U_t(z, y), x) ~~ (\text{commutativity of } U_t)\\
&=U_t(z, U_t(y, x))~~ (\text{by 2.1.3})\\
&=U_t(U_t(x, y), z)~~ (\text{commutativity of } U_t).
\end{align*}

\item[3.1.3.] If $z\in I_{e}$, then $U_{t}(x, U_{t}(y, z))=U_{t}(x, z)$ and
$U_{t}(U_{t}(x, y), z)=U_{t}(x, z)$, implying that $U_t(x, U_t(y, z))=U_t(U_t(x, y),z)$.
\end{itemize}
\item[3.2.] $y\in (e, 1]$.
\begin{itemize}
\item[3.2.1.] If $z\in [0, e)$, then
\begin{align*}
U_t(x, U_t(y, z))&=U_t(U_t(z, y), x)~~ (\text{commutativity of } U_t)\\
&= U_t(z, U_t(y, x))~~ (\text{by 1.2.3})\\
&= U_t(U_t(x, y), z)~~ (\text{commutativity of } U_t).
\end{align*}

\item[3.2.2.] If $z\in (e, 1]$, then
\begin{align*}
U_t(x, U_t(y, z))&= U_t(U_t(z, y), x)~~ (\text{commutativity of } U_t)\\
&= U_t(z, U_t(y, x))~~ (\text{by 2.2.3})\\
&= U_t(U_t(x, y), z)~~ (\text{commutativity of } U_t).
\end{align*}

\item[3.2.3.] If $z\in I_{e}$, then $U_{t}(x, U_{t}(y, z))=U_{t}(x, 1)=1$ and
$U_{t}(U_{t}(x, y), z)=U_{t}(1, z)=1$, implying that $U_t(x, U_t(y, z))=U_t(U_t(x, y),z)$.
\end{itemize}

\item[3.3.] $y\in I_{e}$.
\begin{itemize}
\item[3.3.1.] If $z\in [0, e)$, then
\begin{align*}
U_t(x, U_t(y, z))&=U_t(U_t(z, y), x)~~ (\text{commutativity of } U_t)\\
&= U_t(z, U_t(y, x))~~ (\text{by 1.3.3})\\
&= U_t(U_t(x, y), z)~~ (\text{commutativity of } U_t).
\end{align*}

\item[3.3.2.] If $z\in (e, 1]$, then
\begin{align*}
U_t(x, U_t(y, z))&= U_t(U_t(z, y), x)~~ (\text{commutativity of } U_t)\\
&= U_t(z, U_t(y, x))~~ (\text{by 2.3.3})\\
&= U_t(U_t(x, y), z)~~ (\text{commutativity of } U_t).
\end{align*}

\item[3.3.3.] If $z\in I_{e}$, then
\begin{align*}
U_{t}(x, U_{t}(y, z))&=U_{t}(x, y\vee z)\\
&=\begin{cases}
x\vee y\vee z, & y\vee z\in I_{e}, \\
1, & y\vee z=1,
\end{cases}\\
&=x\vee y \vee z,
\end{align*} and
\begin{align*}
U_{t}(U_{t}(x, y), z)&=U_{t}(x\vee y, z)\\
&=\begin{cases}
x\vee y\vee z, & x\vee y\in I_{e}, \\
1, & x\vee y=1,
\end{cases}\\
&=x\vee y \vee z,
\end{align*}
implying that $U_t(x, U_t(y, z))=U_t(U_t(x, y),z)$.
\end{itemize}
\end{itemize}
\end{itemize}
\end{itemize}
\epf

\begin{theorem}\label{U(s)-operation-Thm}
Let $(X, \leq, 0, 1)$ be a bounded lattice and $e\in L\setminus \{0, 1\}$.
If $S_{e}$ is a $t$-conorm on $[e, 1]$ and the function $U_{s}: L\times L\longrightarrow L$
is defined by
\begin{equation}\label{U(s)-operation*}
U_{s}(x, y)=
\begin{cases}
S_{e}(x, y), & (x, y)\in [e, 1]^2, \\
x\wedge y, & (x, y)\in \left([0, e)\times [e, 1]\right) \cup \left([e, 1]\times [0, e)\right) \cup \left(I_{e}\times I_{e}\right), \\
y, & (x, y)\in [e, 1]\times I_{e}, \\
x, & (x, y)\in I_{e}\times [e, 1], \\
0, & \text{otherwise},
\end{cases}
\end{equation}
then the following statements are equivalent:
\begin{enumerate}[{\rm (i)}]
\item $U_{s}$ is a uninorm on $L$ with the neutral element $e$;

\item $I_{e}=\emptyset$; otherwise, for any $y, z\in I_{e}$, it holds that $y\wedge z\in I_{e}$ or $y\wedge z=0$;

\item $I_{e}\wedge I_{e}\subset I_{e}\cup \{0\}$;

\item $\wedge$ is a $t$-norm on $I_{e} \cup \{0, 1\}$.
\end{enumerate}
%%$I_{e} \vee I_{e} \subset I_{e}$ or $I_{e}\vee I_{e}=\{1\}$.
\end{theorem}
\pf
It can be proved similarly to the proof of Theorem \ref{U(t)-operation-Thm}.
\epf

\section{Characterizations of binary operations defined by \ref{U(T)-operation}}

\begin{theorem}\label{U(T)-operation-Thm}
Let $(X, \leq, 0, 1)$ be a bounded lattice, $e\in L\setminus \{0, 1\}$, and $T_{e}$ be a $t$-norm on $[0, e]$.
Then, the following statements are equivalent:
\begin{enumerate}[{\rm (i)}]
\item\label{1} $U_{T}$ defined by \eqref{U(T)-operation} is a uninorm on $L$ with the neutral
element $e$;

\item\label{2} $I_{e}=\emptyset$; otherwise, for any $y, z\in I_{e}$, it holds that $y\vee z\in I_{e}$ or $y\vee z=1$;

\item\label{3} $I_{e}\vee I_{e}\subset I_{e}\cup \{1\}$;

\item\label{4} $\vee$ is a $t$-conorm on $I_{e} \cup \{0, 1\}$.
\end{enumerate}
%%$I_{e} \vee I_{e} \subset I_{e}$ or $I_{e}\vee I_{e}=\{1\}$.
\end{theorem}
\pf
Clearly, (\ref{2})$\Longleftrightarrow$(\ref{3})$\Longleftrightarrow$(\ref{4}). It remains to show that
(\ref{1})$\Longleftrightarrow$(\ref{2}).

\medskip

(\ref{1})$\Longrightarrow$(\ref{2}).
Suppose, on the contrary, that there exist $y_0, z_0\in I_{e}$ such that $y_0\vee z_0\notin I_{e}$
and $y_0\vee z_0<1$. This implies that $y_{0}\vee z_0\in (e, 1)$. %This implies that
%%% $y_{0}\vee z_0>e$
%\begin{equation}\label{eq-1}
%y_0\vee z_0\in (e, 1).
%\end{equation}
Then,
$$
U_{T}(y_0\vee z_0, U_{T}(y_0, z_0))=U_{T}(y_0\vee z_0, y_0\vee z_0)=1,
$$
and
$$
U_T(U_T(y_0\vee z_0, y_0), z_0)=U_{T}(y_0\vee z_0, z_0)=y_0\vee z_0<1,
$$
which contradict with the associativity of $U_{T}$.

\medskip

(\ref{2})$\Longrightarrow$(\ref{1}).
Clearly, $U_{T}$ is a uninorm when $I_{e}=\emptyset$. Without loss of generality, assume
that $I_{e}\neq \emptyset$. Firstly, it is not difficult to verify
that $U_T$ is a commutative binary operation with the neutral element $e$.

i.1) {\it Monotonicity}. For any $x, y\in L$ with $x\leq y$, it holds that $U_T(x, z)\leq U_T(y, z)$ for all $z\in L$.
Consider the following cases:
\begin{itemize}
\item[1.] $x\in [0, e)$.
\begin{itemize}
\item[1.1.] If $y\in [0, e)$, then
\begin{equation}\label{eq-2}
U_{T}(x, z)=
\begin{cases}
T_{e}(x, z), & z\in [0, e), \\
x, & z=e, \\
z, & z\in (e, 1], \\
z, & z\in I_{e},
\end{cases}
\end{equation}
and
$$
U_{T}(y, z)=
\begin{cases}
T_{e}(y, z), & z\in [0, e), \\
y, & z=e, \\
z, & z\in (e, 1], \\
z, & z\in I_{e},
\end{cases}
$$
implying that $U_{T}(x, z)\leq U_{T}(y, z)$, because $T_{e}$ is increasing.

\item[1.2.] If $y\in (e, 1]$, then,
$$
U_{T}(y, z)=
\begin{cases}
y, & z\in [0, e), \\
y, & z=e, \\
1, & z\in (e, 1], \\
y\vee z, & z\in I_{e}.
\end{cases}
$$
This, together with \eqref{eq-2} and $T_{e}(x, z)\leq x$, implies that $U_T(x, z)\leq U_T(y, z)$.

\item[1.3.] If $y\in I_{e}$, then
\begin{equation}\label{eq-5}
U_{T}(y, z)=
\begin{cases}
y, & z\in [0, e), \\
y, & z=e, \\
y\vee z, & z\in (e, 1], \\
y\vee z, & z\in I_{e}.
\end{cases}
\end{equation}
This, together with \eqref{eq-2} and $T_{e}(x, z)\leq x$, implies that $U_T(x, z)\leq U_T(y, z)$.

\item[1.4.] If $y=e$, then $U_{T}(y, z)=z$. This, together with \eqref{eq-2}, implies that $U_T(x, z)\leq U_T(y, z)$.
\end{itemize}

\item[2.] $x\in (e, 1]$. From $x\leq y$, it follows that $y\in (e, 1]$. Then,
\begin{equation}\label{eq-3}
U_T(x, z)=
\begin{cases}
x, & z\in [0, e), \\
x, & z=e, \\
1, & z\in (e, 1], \\
x\vee z, & z\in I_{e},
\end{cases}
\end{equation}
and
\begin{equation}\label{eq-4}
U_T(y, z)=
\begin{cases}
y, & z\in [0, e), \\
y, & z=e, \\
1, & z\in (e, 1], \\
y\vee z, & z\in I_{e},
\end{cases}
\end{equation}
implying that $U_T(x, z)\leq U_T(y, z)$.

\item[3.] $x=e$. This implies that $y\in [e, 1]$. Clearly, $U_T(x, z)=z= U_T(y, z)$ if $y=e$. If $y\in (e, 1]$,
applying \eqref{eq-4} and $U_{T}(x, z)=z$, it follows that $U_{T}(x, z)\leq U_{T}(y, z)$.

\item[4.] $x\in I_{e}$. From $x\leq y$, it follows that $y\in (e, 1] \cup I_{e}$. Applying \eqref{U(T)-operation}
yields that
\begin{equation}\label{1-16}
U_{T}(x, z)=
\begin{cases}
x, & z\in [0, e), \\
x, & z=e, \\
x\vee z, & z\in (e, 1], \\
x\vee z, & z\in I_{e}.
\end{cases}
\end{equation}
\begin{itemize}
\item[4.1.] If $y\in (e, 1]$, from \eqref{eq-4} and \eqref{1-16}, it follows that $U_{T}(x, z)\leq U_{T}(y, z)$.

\item[4.2.] If $y\in I_{e}$, from \eqref{eq-5}  and \eqref{1-16}, it is clear that $U_{T}(x, z)\leq U_{T}(y, z)$.
\end{itemize}
\end{itemize}

i.2) {\it Associativity.} For any $x, y, z\in L$, it holds that $U_T(x, U_T(y, z))=U_T(U_T(x, y),z)$.%% for any $x, y, z\in L$.

If one of the elements $x$, $y$ and $z$ is equal to $e$, then the equality is always satisfied.
Otherwise, consider the following cases:
\begin{itemize}
\item[1.] $x\in [0, e)$.
\begin{itemize}
\item[1.1.] $y\in [0, e)$.
\begin{itemize}
\item[1.1.1.] If $z\in [0, e)$, then $U_{T}(x, U_{T}(y, z))=U_{T}(x, T_e(y, z))=T_{e}(x, T_{e}(y, z))=
T_{e}(T_{e}(x, y), z)=T_{e}(U_{T}(x, y), z)=U_{T}(U_{T}(x, y), z)$.

\item[1.1.2.] If $z\in (e, 1]$, then $U_{T}(x, U_{T}(y, z))=U_{T}(x, z)=z$ and $U_{T}(U_{T}(x, y), z)
=U_{T}(T_{e}(x, y), z)=z$, implying that $U_T(x, U_T(y, z))=U_T(U_T(x, y),z)$.

\item[1.1.3.] If $z\in I_{e}$, then $U_{T}(x, U_{T}(y, z))=U_{T}(x, z)=z$ and
$U_{T}(U_{T}(x, y), z)=U_{T}(T_{e}(x, y), z)=z$, implying that $U_T(x, U_T(y, z))=U_T(U_T(x, y),z)$.
\end{itemize}
\item[1.2.] $y\in (e, 1]$.
\begin{itemize}
\item[1.2.1.] If $z\in [0, e)$, then $U_{T}(x, U_{T}(y, z))=U_{T}(x, y)=y$ and
$U_{T}(U_{T}(x, y), z)=U_{T}(y, z)=y$, implying that $U_T(x, U_T(y, z))=U_T(U_T(x, y),z)$.

\item[1.2.2.] If $z\in (e, 1]$, then $U_{T}(x, U_{T}(y, z))=U_{T}(x, 1)=1$ and
$U_{T}(U_{T}(x, y), z)=U_{T}(y, z)=1$, implying that $U_T(x, U_T(y, z))=U_T(U_T(x, y),z)$.

\item[1.2.3.] If $z\in I_{e}$, then $U_{T}(x, U_{T}(y, z))=U_{T}(x, y\vee z)=y\vee z$ and
$U_{T}(U_{T}(x, y), z)=U_{T}(y, z)=y\vee z$, implying that $U_T(x, U_T(y, z))=U_T(U_T(x, y),z)$.
\end{itemize}
\item[1.3.] $y\in I_{e}$.
\begin{itemize}
\item[1.3.1.] If $z\in [0, e)$, then $U_{T}(x, U_{T}(y, z))=U_{T}(x, y)=y$ and
$U_{T}(U_{T}(x, y), z)=U_{T}(y, z)=y$, implying that $U_T(x, U_T(y, z))=U_T(U_T(x, y),z)$.

\item[1.3.2.] If $z\in (e, 1]$, then $U_{T}(x, U_{T}(y, z))=U_{T}(x, y\vee z)=y\vee z$ and
$U_{T}(U_{T}(x, y), z)=U_{T}(y, z)=y\vee z$, implying that $U_T(x, U_T(y, z))=U_T(U_T(x, y),z)$.

\item[1.3.3.] If $z\in I_{e}$, then
\begin{align*}
U_{T}(x, U_{T}(y, z))
&=U_{T}(x, y\vee z)\\
&=
\begin{cases}
y\vee z, & y\vee z\in I_{e}, \\
1, & y\vee z=1,
\end{cases}\\
&=y\vee z,
\end{align*}
and
$U_{T}(U_{T}(x, y), z)=U_{T}(y, z)=y\vee z$, implying that $U_T(x, U_T(y, z))=U_T(U_T(x, y),z)$.
\end{itemize}
\end{itemize}

\item[2.] $x\in (e, 1]$.
\begin{itemize}
\item[2.1.] $y\in [0, e)$.
\begin{itemize}
\item[2.1.1.] If $z\in [0, e)$, then
\begin{align*}
U_T(x, U_T(y, z))&=U_T(U_T(z, y), x) ~~ (\text{commutativity of } U_T)\\
&=U_T(z, U_T(y, x))~~ (\text{by 1.1.2})\\
&=U_T(U_T(x, y), z)~~ (\text{commutativity of } U_T).
\end{align*}

\item[2.1.2.] If $z\in (e, 1]$, then $U_{T}(x, U_{T}(y, z))=U_{T}(x, z)$ and
$U_{T}(U_{T}(x, y), z)=U_{T}(x, z)$, implying that $U_T(x, U_T(y, z))=U_T(U_T(x, y),z)$.

\item[2.1.3.] If $z\in I_{e}$, then $U_{T}(x, U_{T}(y, z))=U_{T}(x, z)$ and
$U_{T}(U_{T}(x, y), z)=U_{T}(x, z)$, implying that $U_T(x, U_T(y, z))=U_T(U_T(x, y),z)$.
\end{itemize}
\item[2.2.] $y\in (e, 1]$.
\begin{itemize}
\item[2.2.1.] If $z\in [0, e)$, then
\begin{align*}
U_T(x, U_T(y, z))&=U_T(U_T(z, y), x)~~ (\text{commutativity of } U_T)\\
&= U_T(z, U_T(y, x))~~ (\text{by 1.2.2})\\
&= U_T(U_T(x, y), z)~~ (\text{commutativity of } U_T).
\end{align*}

\item[2.2.2.] If $z\in (e, 1]$, then $U_{T}(x, U_{T}(y, z))=U_{T}(x, 1)=1$ and
$U_{T}(U_{T}(x, y), z)=U_{T}(1, z)=1$, implying that $U_T(x, U_T(y, z))=U_T(U_T(x, y),z)$.

\item[2.2.3.] If $z\in I_{e}$, then $U_{T}(x, U_{T}(y, z))=U_{T}(x, y\vee z)=1$ and
$U_{T}(U_{T}(x, y), z)=U_{T}(1, z)=1$, implying that $U_T(x, U_T(y, z))=U_T(U_T(x, y),z)$.
\end{itemize}

\item[2.3.] $y\in I_{e}$.
\begin{itemize}
\item[2.3.1.] If $z\in [0, e)$, then \begin{align*}
U_T(x, U_T(y, z))&=U_T(U_T(z, y), x) ~~ (\text{commutativity of } U_T)\\
&=U_T(z, U_T(y, x))~~ (\text{by 1.3.2})\\
&=U_T(U_T(x, y), z)~~ (\text{commutativity of } U_T).
\end{align*}

\item[2.3.2.] If $z\in (e, 1]$, then $U_{T}(x, U_{T}(y, z))=U_{T}(x, y\vee z)=1$ and
$U_{T}(U_{T}(x, y), z)=U_{T}(x\vee y, z)=1$, implying that $U_T(x, U_T(y, z))=U_T(U_T(x, y),z)$.

\item[2.3.3.] If $z\in I_{e}$, then
\begin{align*}
U_{T}(x, U_{T}(y, z))&=U_{T}(x, y\vee z)\\
&=\begin{cases}
x\vee y\vee z, & y\vee z\in I_{e}, \\
1, & y\vee z=1,
\end{cases}\\
& =x\vee y\vee z,
\end{align*}
and
$U_{T}(U_{T}(x, y), z)=U_{T}(x\vee y, z)=x\vee y\vee z$,
implying that $U_T(x, U_T(y, z))=U_T(U_T(x, y),z)$.
\end{itemize}

\item[3.] $x\in I_{e}$.
\begin{itemize}
\item[3.1.] $y\in [0, e)$.
\begin{itemize}
\item[3.1.1.] If $z\in [0, e)$, then
\begin{align*}
U_T(x, U_T(y, z))&=U_T(U_T(z, y), x) ~~ (\text{commutativity of } U_T)\\
&=U_T(z, U_T(y, x))~~ (\text{by 1.1.3})\\
&=U_T(U_T(x, y), z)~~ (\text{commutativity of } U_T).
\end{align*}

\item[3.1.2.] If $z\in (e, 1]$, then
\begin{align*}
U_T(x, U_T(y, z))&=U_T(U_T(z, y), x) ~~ (\text{commutativity of } U_T)\\
&=U_T(z, U_T(y, x))~~ (\text{by 2.1.3})\\
&=U_T(U_T(x, y), z)~~ (\text{commutativity of } U_T).
\end{align*}

\item[3.1.3.] If $z\in I_{e}$, then $U_{T}(x, U_{T}(y, z))=U_{T}(x, z)$ and
$U_{T}(U_{T}(x, y), z)=U_{T}(x, z)$, implying that $U_T(x, U_T(y, z))=U_T(U_T(x, y),z)$.
\end{itemize}
\item[3.2.] $y\in (e, 1]$.
\begin{itemize}
\item[3.2.1.] If $z\in [0, e)$, then
\begin{align*}
U_T(x, U_T(y, z))&=U_T(U_T(z, y), x)~~ (\text{commutativity of } U_T)\\
&= U_T(z, U_T(y, x))~~ (\text{by 1.2.3})\\
&= U_T(U_T(x, y), z)~~ (\text{commutativity of } U_T).
\end{align*}

\item[3.2.2.] If $z\in (e, 1]$, then
\begin{align*}
U_T(x, U_T(y, z))&= U_T(U_T(z, y), x)~~ (\text{commutativity of } U_T)\\
&= U_T(z, U_T(y, x))~~ (\text{by 2.2.3})\\
&= U_T(U_T(x, y), z)~~ (\text{commutativity of } U_T).
\end{align*}

\item[3.2.3.] If $z\in I_{e}$, then $U_{T}(x, U_{T}(y, z))=U_{T}(x, y\vee z)=x\vee y\vee z$ and
$U_{T}(U_{T}(x, y), z)=U_{T}(x\vee y, z)=x\vee y\vee z$, implying that $U_T(x, U_T(y, z))=U_T(U_T(x, y),z)$.
\end{itemize}

\item[3.3.] $y\in I_{e}$.
\begin{itemize}
\item[3.3.1.] If $z\in [0, e)$, then
\begin{align*}
U_T(x, U_T(y, z))&=U_T(U_T(z, y), x)~~ (\text{commutativity of } U_T)\\
&= U_T(z, U_T(y, x))~~ (\text{by 1.3.3})\\
&= U_T(U_T(x, y), z)~~ (\text{commutativity of } U_T).
\end{align*}

\item[3.3.2.] If $z\in (e, 1]$, then
\begin{align*}
U_T(x, U_T(y, z))&= U_T(U_T(z, y), x)~~ (\text{commutativity of } U_T)\\
&= U_T(z, U_T(y, x))~~ (\text{by 2.3.3})\\
&= U_T(U_T(x, y), z)~~ (\text{commutativity of } U_T).
\end{align*}

\item[3.3.3.] If $z\in I_{e}$, then
\begin{align*}
U_{T}(x, U_{T}(y, z))&=U_{T}(x, y\vee z)\\
&=\begin{cases}
x\vee y\vee z, & y\vee z\in I_{e}, \\
1, & y\vee z=1,
\end{cases}\\
&=x\vee y \vee z,
\end{align*}
and
\begin{align*}
U_{T}(U_{T}(x, y), z)&=U_{T}(x\vee y, z)\\
&=\begin{cases}
x\vee y\vee z, & x\vee y\in I_{e}, \\
1, & x\vee y=1,
\end{cases}\\
&=x\vee y \vee z,
\end{align*}
implying that $U_T(x, U_T(y, z))=U_T(U_T(x, y),z)$.
\end{itemize}
\end{itemize}
\end{itemize}
\end{itemize}
\epf

%\begin{theorem}\label{U(S)-operation-Thm}
%Let $(X, \leq, 0, 1)$ be a bounded lattice, $e\in L\setminus \{0, 1\}$, and $S_{e}$ be a $t$-conorm on $[e, 1]$.
%Then, the following statements are equivalent:
%\begin{enumerate}[{\rm (i)}]
%\item $U_{S}$ defined by \eqref{U(S)-operation} is a uninorm on $L$ with the neutral
%element $e$;
%
%\item $I_{e}=\emptyset$; otherwise, for any $y, z\in I_{e}$, $y\wedge z\in I_{e}$ or $y\wedge z=0$;
%
%\item $I_{e}\wedge I_{e}\subset I_{e}\cup \{0\}$;
%
%\item $\wedge$ is a $t$-norm on $I_{e} \cup \{0, 1\}$.
%\end{enumerate}
%%%$I_{e} \vee I_{e} \subset I_{e}$ or $I_{e}\vee I_{e}=\{1\}$.
%\end{theorem}
%\pf
%It can be proved as dual of Theorem \ref{U(T)-operation-Thm}.
%\epf

\section{Characterizations of uninorms constructed by A{\c{s}}{\i}c{\i} and Mesiar \cite{AM2020}}
This section is devoted to characterizing the uninorm property for the operations introduced by
A{\c{s}}{\i}c{\i} and Mesiar \cite{AM2020} (see Theorems~\ref{U(Te)-operation-Thm} and
\ref{U(Se)-operation-Thm}). First, the incomparability is defined for
two subsets of a lattice.
\begin{definition}
Let $(L, \leq, 0, 1)$ be a bounded lattice and $A, B\subset L$.
Then, $A$ and $B$ are {\it comparable} if there exist $x\in A$ and $y\in B$
such that $x$ and $y$ are comparable. Otherwise, $A$ and $B$ are {\it incomparable},
i.e., $x\| y$ for all $x\in A$ and $y\in B$, and in this case, it is denoted as $A\| B$.
\end{definition}

Recently, A{\c{s}}{\i}c{\i} and Mesiar \cite{AM2020} obtained the following results:
\begin{theorem}{\rm \cite[Theorem~6]{AM2020}}\label{AM-Thm-6}
Let $(L, \leq, 0, 1)$ be a bounded lattice with $e\in L\setminus \{0, 1\}$ such that for all
$x\in I_e$ and $y\in (0, e]$ it holds $x\| y$, i.e., $I_{e}\| (0, e]$. If $T_{e}$ is a $t$-norm on $[0, e]$,
then the function $U_{T_e}: L^2\longrightarrow L$ defined by
\begin{equation}\label{U(Te)-operation}
U_{T_{e}}(x, y)=
\begin{cases}
T_{e}(x, y), & (x, y)\in [0, e]^2, \\
y, & (x, y)\in [e, 1]\times I_{e}, \\
x, & (x, y)\in I_{e}\times [e, 1], \\
0, & (x, y)\in \left([0, e)\times I_{e}\right) \cup \left(I_{e}\times [0, e)\right), \\
x\wedge y, & (x, y)\in \left([0, e)\times [e, 1]\right) \cup \left([e, 1]\times [0, e)\right) \cup \left(I_{e}\times I_{e}\right), \\
x\vee y, & \text{otherwise},
\end{cases}
\end{equation}
is a uninorm on $L$ with the neutral element $e$.
\end{theorem}

\begin{theorem}{\rm \cite[Theorem~7]{AM2020}}\label{AM-Thm-7}
Let $(L, \leq, 0, 1)$ be a bounded lattice with $e\in L\setminus \{0, 1\}$ such that for all
$x\in I_e$ and $y\in (0, e]$ it holds $x\| y$, i.e., $I_{e}\| (0, e]$. If $S_{e}$ is a $t$-conorm on $[e, 1]$,
then the function $U_{S_e}: L^2\longrightarrow L$ defined by
\begin{equation}\label{U(Se)-operation}
U_{S_{e}}(x, y)=
\begin{cases}
S_{e}(x, y), & (x, y)\in [e, 1]^2, \\
y, & (x, y)\in [e, 1]\times I_{e}, \\
x, & (x, y)\in I_{e}\times [e, 1], \\
0, & (x, y)\in \left([0, e)\times I_{e}\right) \cup \left(I_{e}\times [0, e)\right), \\
x\wedge y, & (x, y)\in \left([0, e)\times [e, 1]\right) \cup \left([e, 1]\times [0, e)\right) \cup \left(I_{e}\times I_{e}\right), \\
x\vee y, & \text{otherwise},
\end{cases}
\end{equation}
is a uninorm on $L$ with the neutral element $e$.
\end{theorem}

\begin{theorem}\label{U(Te)-operation-Thm}
Let $(X, \leq, 0, 1)$ be a bounded lattice, $e\in L\setminus \{0, 1\}$, and
$T_{e}$ be a $t$-norm on $[0, e]$.
Then, the following statements are equivalent:
\begin{enumerate}[{\rm (i)}]
\item\label{1.i} $U_{T_{e}}$ defined by \eqref{U(Te)-operation} is a uninorm on $L$ with the neutral
element $e$;
%
%\item\label{1.i'} $U_{S_{e}}$ defined by \eqref{U(Se)-operation} is a uninorm on $L$ with the neutral
%element $e$;

\item\label{1.ii} The following hold:
\begin{itemize}
\item[{\rm ii-1)}] $\mathscr{P}=\emptyset$ or $T_{e}|_{\mathscr{P}\times [0, e) \cup [0, e)\times \mathscr{P}}\equiv 0$,
where $\mathscr{P}=\{x\in (0, e): \exists y\in I_{e}, ~ x\leq y\}$;

\item[{\rm ii-2)}] $I_{e}\wedge I_{e}\subset I_{e}\cup \{0\}$.
\end{itemize}
%\item\label{1.iii} $I_{e}\wedge I_{e}\subset I_{e}\cup \{0\}$;
%
%\item\label{1.iv} $\wedge$ is a $t$-norm on $I_{e} \cup \{0, 1\}$.
\end{enumerate}
%%$I_{e} \vee I_{e} \subset I_{e}$ or $I_{e}\vee I_{e}=\{1\}$.
\end{theorem}
\pf
(\ref{1.i})$\Longrightarrow$(\ref{1.ii}).

(i)$\Longrightarrow$ii.1). Suppose, on the contrary, that there exist $x_0\in [0, e)$ and $y_0\in \mathscr{P}$ such that
$T_{e}(x_0, y_0)>0$. From $y_0\in \mathscr{P}$, it follows that there exists $z\in I_{e}$
such that $y_{0}\leq z_0$. Since $U_{T_{e}}$ is increasing, applying \eqref{U(Te)-operation} yields that
$$
0<T_{e}(x_{0}, y_0)=U_{T_{e}}(x_0, y_0)\leq U_{T_{e}}(x_0, z_0)=0,
$$
which is a contradiction. %% with the associativity of $U_{T_{e}}$.

(i)$\Longrightarrow$ii.2). Suppose that $I_{e}\wedge I_{e}\nsubseteq I_{e}\cup \{0\}$, i.e.,
there exist $y_1, z_1\in I_{e}$ such that $0< y_1\wedge z_1$ and $y_1\wedge z_1\notin I_{e}$. This implies that
$y_1\wedge z_1\in (0, e)$. From~\eqref{U(Te)-operation} and the associativity of $U_{T_{e}}$, it follows that
$$
0=U_{T_{e}}(y_1, U_{T_{e}}(y_1, z_1))=U_{T_{e}}( U_{T_{e}}(y_1,y_1), z_1)=y_1\wedge z_1>0,
$$
which is a contradiction.

\medskip

(\ref{1.ii})$\Longrightarrow$(\ref{1.i}).
When $\mathscr{P}=\emptyset$, this has been proved in \cite[Theorem~6]{AM2020}. Now assume that
$T_{e}|_{\mathscr{P}\times [0, e) \cup [0, e)\times \mathscr{P}}\equiv 0$ and $I_{e}\wedge I_{e}\subset I_{e}\cup \{0\}$.
Firstly, it is not difficult to verify that $U_{T_e}$ is a commutative binary operation with the neutral element $e$.

i.1) {\it Monotonicity}. For any $x, y\in L$ with $x\leq y$, it holds that $U_{T_e}(x, z)\leq U_{T_e}(y, z)$ for all $z\in L$.
Consider the following cases:
\begin{itemize}
\item[1.] $x\in [0, e)$.
\begin{itemize}
\item[1.1.] If $y\in [0, e)$, then
\begin{equation}\label{eq-2**}
U_{T_e}(x, z)=
\begin{cases}
T_{e}(x, z), & z\in [0, e), \\
x, & z=e, \\
x, & z\in (e, 1], \\
0, & z\in I_{e},
\end{cases}
\end{equation}
and
$$
U_{T_e}(y, z)=
\begin{cases}
T_{e}(y, z), & z\in [0, e), \\
y, & z=e, \\
y, & z\in (e, 1], \\
0, & z\in I_{e},
\end{cases}
$$
implying that $U_{T_e}(x, z)\leq U_{T_e}(y, z)$, because $T_{e}$ is increasing.

\item[1.2.] If $y\in (e, 1]$, then,
$$
U_{T_e}(y, z)=
\begin{cases}
z, & z\in [0, e), \\
y, & z=e, \\
y\vee z, & z\in (e, 1], \\
z, & z\in I_{e}.
\end{cases}
$$
This, together with \eqref{eq-2**} and $T_{e}(x, z)\leq z$, implies that $U_{T_e}(x, z)\leq U_{T_e}(y, z)$.

\item[1.3.] If $y\in I_{e}$, then $x\in \mathscr{P}$. Applying \eqref{eq-2**} and \eqref{U(Te)-operation}
follows that
$$
U_{T_e}(x, z)=
\begin{cases}
0, & z\in [0, e), \\
x, & z=e, \\
x, & z\in (e, 1], \\
0, & z\in I_{e},
\end{cases}
$$
and
\begin{equation}\label{eq-5**}
U_{T_e}(y, z)=
\begin{cases}
0, & z\in [0, e), \\
y, & z=e, \\
y, & z\in (e, 1], \\
y\wedge z, & z\in I_{e},
\end{cases}
\end{equation}
implying that $U_{T_e}(x, z)\leq U_{T_e}(y, z)$.

\item[1.4.] If $y=e$, then $U_{T_e}(y, z)=z$. This, together with \eqref{eq-2**}, implies that $U_{T_e}(x, z)\leq U_{T_e}(y, z)$.
\end{itemize}

\item[2.] $x\in (e, 1]$. From $x\leq y$, it follows that $y\in (e, 1]$. Then,
\begin{equation}\label{eq-3**}
U_{T_e}(x, z)=
\begin{cases}
z, & z\in [0, e), \\
x, & z=e, \\
x\vee z, & z\in (e, 1], \\
z, & z\in I_{e},
\end{cases}
\end{equation}
and
\begin{equation}\label{eq-4**}
U_{T_e}(y, z)=
\begin{cases}
z, & z\in [0, e), \\
y, & z=e, \\
y\vee z, & z\in (e, 1], \\
z, & z\in I_{e},
\end{cases}
\end{equation}
implying that $U_{T_e}(x, z)\leq U_{T_e}(y, z)$.

\item[3.] $x=e$. This implies that $y\in [e, 1]$. Clearly, $U_{T_e}(x, z)=z= U_{T_e}(y, z)$ if $y=e$. If $y\in (e, 1]$,
applying \eqref{eq-4**} and $U_{T_e}(x, z)=z$, it follows that $U_{T_e}(x, z)\leq U_{T_e}(y, z)$.

\item[4.] $x\in I_{e}$. From $x\leq y$, it follows that $y\in (e, 1] \cup I_{e}$. Applying \eqref{U(Te)-operation}
yields that
\begin{equation}\label{1-16**}
U_{T_e}(x, z)=
\begin{cases}
0, & z\in [0, e), \\
x, & z=e, \\
x, & z\in (e, 1], \\
x\wedge z, & z\in I_{e}.
\end{cases}
\end{equation}
\begin{itemize}
\item[4.1.] If $y\in (e, 1]$, from \eqref{eq-4**} and \eqref{1-16**}, it follows that $U_{T_e}(x, z)\leq U_{T_e}(y, z)$.

\item[4.2.] If $y\in I_{e}$, from \eqref{eq-5**}  and \eqref{1-16**}, it is clear that $U_{T_e}(x, z)\leq U_{T_e}(y, z)$.
\end{itemize}
\end{itemize}

i.2) {\it Associativity.} For any $x, y, z\in L$, it holds that $U_{T_e}(x, U_{T_e}(y, z))=U_{T_e}(U_{T_e}(x, y),z)$.%% for any $x, y, z\in L$.

If one of the elements $x$, $y$ and $z$ is equal to $e$, then the equality is always satisfied.
Otherwise, consider the following cases:
\begin{itemize}
\item[1.] $x\in [0, e)$.
\begin{itemize}
\item[1.1.] $y\in [0, e)$.
\begin{itemize}
\item[1.1.1.] If $z\in [0, e)$, then $U_{T_e}(x, U_{T_e}(y, z))=U_{T_e}(x, T_e(y, z))=T_{e}(x, T_{e}(y, z))=
T_{e}(T_{e}(x, y), z)=T_{e}(U_{T_e}(x, y), z)=U_{T_e}(U_{T_e}(x, y), z)$.

\item[1.1.2.] If $z\in (e, 1]$, then $U_{T_e}(x, U_{T_e}(y, z))=U_{T_e}(x, y)=T_{e}(x, y)$ and $U_{T_e}(U_{T_e}(x, y), z)
=U_{T_e}(T_{e}(x, y), z)=T_{e}(x, y)$, implying that $U_{T_e}(x, U_{T_e}(y, z))=U_{T_e}(U_{T_e}(x, y),z)$.

\item[1.1.3.] If $z\in I_{e}$, then $U_{T_e}(x, U_{T_e}(y, z))=U_{T_e}(x, 0)=0$ and
$U_{T_e}(U_{T_e}(x, y), z)=U_{T_e}(T_{e}(x, y), z)=0$, implying that $U_{T_e}(x, U_{T_e}(y, z))=U_{T_e}(U_{T_e}(x, y),z)$.
\end{itemize}
\item[1.2.] $y\in (e, 1]$.
\begin{itemize}
\item[1.2.1.] If $z\in [0, e)$, then $U_{T_e}(x, U_{T_e}(y, z))=U_{T_e}(x, z)=T_{e}(x, z)$ and
$U_{T_e}(U_{T_e}(x, y), z)=U_{T_e}(x, z)=T_{e}(x, z)$, implying that $U_{T_e}(x, U_{T_e}(y, z))=U_{T_e}(U_{T_e}(x, y),z)$.

\item[1.2.2.] If $z\in (e, 1]$, then $U_{T_e}(x, U_{T_e}(y, z))=U_{T_e}(x, y\vee z)=x$ and
$U_{T_e}(U_{T_e}(x, y), z)=U_{T_e}(x, z)=x$, implying that $U_{T_e}(x, U_{T_e}(y, z))=U_{T_e}(U_{T_e}(x, y),z)$.

\item[1.2.3.] If $z\in I_{e}$, then $U_{T_e}(x, U_{T_e}(y, z))=U_{T_e}(x, z)=0$ and
$U_{T_e}(U_{T_e}(x, y), z)=U_{T_e}(x, z)=0$, implying that $U_{T_e}(x, U_{T_e}(y, z))=U_{T_e}(U_{T_e}(x, y),z)$.
\end{itemize}
\item[1.3.] $y\in I_{e}$.
\begin{itemize}
\item[1.3.1.] If $z\in [0, e)$, then $U_{T_e}(x, U_{T_e}(y, z))=U_{T_e}(x, 0)=0$ and
$U_{T_e}(U_{T_e}(x, y), z)=U_{T_e}(0, z)=0$, implying that $U_{T_e}(x, U_{T_e}(y, z))=U_{T_e}(U_{T_e}(x, y),z)$.

\item[1.3.2.] If $z\in (e, 1]$, then $U_{T_e}(x, U_{T_e}(y, z))=U_{T_e}(x, y)=0$ and
$U_{T_e}(U_{T_e}(x, y), z)=U_{T_e}(0, z)=0$, implying that $U_{T_e}(x, U_{T_e}(y, z))=U_{T_e}(U_{T_e}(x, y),z)$.

\item[1.3.3.] If $z\in I_{e}$, then
\begin{align*}
U_{T_e}(x, U_{T_e}(y, z))
&=U_{T_e}(x, y\wedge z)\\
&=
\begin{cases}
0, & y\wedge z\in I_{e}, \\
0, & y\wedge z=0,
\end{cases}\\
&=0,
\end{align*}
and
$U_{T_e}(U_{T_e}(x, y), z)=U_{T_e}(0, z)=0$, implying that $U_{T_e}(x, U_{T_e}(y, z))=U_{T_e}(U_{T_e}(x, y),z)$.
\end{itemize}
\end{itemize}

\item[2.] $x\in (e, 1]$.
\begin{itemize}
\item[2.1.] $y\in [0, e)$.
\begin{itemize}
\item[2.1.1.] If $z\in [0, e)$, then
\begin{align*}
U_{T_e}(x, U_{T_e}(y, z))&=U_{T_e}(U_{T_e}(z, y), x) ~~ (\text{commutativity of } U_{T_e})\\
&=U_{T_e}(z, U_{T_e}(y, x))~~ (\text{by 1.1.2})\\
&=U_{T_e}(U_{T_e}(x, y), z)~~ (\text{commutativity of } U_{T_e}).
\end{align*}

\item[2.1.2.] If $z\in (e, 1]$, then $U_{T_e}(x, U_{T_e}(y, z))=U_{T_e}(x, y)=y$ and
$U_{T_e}(U_{T_e}(x, y), z)=U_{T_e}(y, z)=y$, implying that $U_{T_e}(x, U_{T_e}(y, z))=U_{T_e}(U_{T_e}(x, y),z)$.

\item[2.1.3.] If $z\in I_{e}$, then $U_{T_e}(x, U_{T_e}(y, z))=U_{T_e}(x, 0)=0$ and
$U_{T_e}(U_{T_e}(x, y), z)=U_{T_e}(y, z)=0$, implying that $U_{T_e}(x, U_{T_e}(y, z))=U_{T_e}(U_{T_e}(x, y),z)$.
\end{itemize}
\item[2.2.] $y\in (e, 1]$.
\begin{itemize}
\item[2.2.1.] If $z\in [0, e)$, then
\begin{align*}
U_{T_e}(x, U_{T_e}(y, z))&=U_{T_e}(U_{T_e}(z, y), x)~~ (\text{commutativity of } U_{T_e})\\
&= U_{T_e}(z, U_{T_e}(y, x))~~ (\text{by 1.2.2})\\
&= U_{T_e}(U_{T_e}(x, y), z)~~ (\text{commutativity of } U_{T_e}).
\end{align*}

\item[2.2.2.] If $z\in (e, 1]$, then $U_{T_e}(x, U_{T_e}(y, z))=U_{T_e}(x, y\vee z)=x\vee y\vee z$ and
$U_{T_e}(U_{T_e}(x, y), z)=U_{T_e}(x\vee y, z)=x\vee y\vee z$, implying that $U_{T_e}(x, U_{T_e}(y, z))=U_{T_e}(U_{T_e}(x, y),z)$.

\item[2.2.3.] If $z\in I_{e}$, then $U_{T_e}(x, U_{T_e}(y, z))=U_{T_e}(x, z)=z$ and
$U_{T_e}(U_{T_e}(x, y), z)=U_{T_e}(x\vee y, z)=z$, implying that $U_{T_e}(x, U_{T_e}(y, z))=U_{T_e}(U_{T_e}(x, y),z)$.
\end{itemize}

\item[2.3.] $y\in I_{e}$.
\begin{itemize}
\item[2.3.1.] If $z\in [0, e)$, then \begin{align*}
U_{T_e}(x, U_{T_e}(y, z))&=U_{T_e}(U_{T_e}(z, y), x) ~~ (\text{commutativity of } U_{T_e})\\
&=U_{T_e}(z, U_{T_e}(y, x))~~ (\text{by 1.3.2})\\
&=U_{T_e}(U_{T_e}(x, y), z)~~ (\text{commutativity of } U_{T_e}).
\end{align*}

\item[2.3.2.] If $z\in (e, 1]$, then $U_{T_e}(x, U_{T_e}(y, z))=U_{T_e}(x, y)=y$ and
$U_{T_e}(U_{T_e}(x, y), z)=U_{T_e}(y, z)=y$, implying that $U_{T_e}(x, U_{T_e}(y, z))=U_{T_e}(U_{T_e}(x, y),z)$.

\item[2.3.3.] If $z\in I_{e}$, then
\begin{align*}
U_{T_e}(x, U_{T_e}(y, z))&=U_{T_e}(x, y\wedge z)\\
&=\begin{cases}
y\wedge z, & y\wedge z\in I_{e}, \\
0, & y\wedge z=0,
\end{cases}\\
& =y\wedge z,
\end{align*}
and
$U_{T_e}(U_{T_e}(x, y), z)=U_{T_e}(y, z)=y\wedge z$,
implying that $U_{T_e}(x, U_{T_e}(y, z))=U_{T_e}(U_{T_e}(x, y),z)$.
\end{itemize}

\item[3.] $x\in I_{e}$.
\begin{itemize}
\item[3.1.] $y\in [0, e)$.
\begin{itemize}
\item[3.1.1.] If $z\in [0, e)$, then
\begin{align*}
U_{T_e}(x, U_{T_e}(y, z))&=U_{T_e}(U_{T_e}(z, y), x) ~~ (\text{commutativity of } U_{T_e})\\
&=U_{T_e}(z, U_{T_e}(y, x))~~ (\text{by 1.1.3})\\
&=U_{T_e}(U_{T_e}(x, y), z)~~ (\text{commutativity of } U_{T_e}).
\end{align*}

\item[3.1.2.] If $z\in (e, 1]$, then
\begin{align*}
U_{T_e}(x, U_{T_e}(y, z))&=U_{T_e}(U_{T_e}(z, y), x) ~~ (\text{commutativity of } U_{T_e})\\
&=U_{T_e}(z, U_{T_e}(y, x))~~ (\text{by 2.1.3})\\
&=U_{T_e}(U_{T_e}(x, y), z)~~ (\text{commutativity of } U_{T_e}).
\end{align*}

\item[3.1.3.] If $z\in I_{e}$, then $U_{T_e}(x, U_{T_e}(y, z))=U_{T_e}(x, 0)=0$ and
$U_{T_e}(U_{T_e}(x, y), z)=U_{T_e}(0, z)=0$, implying that $U_{T_e}(x, U_{T_e}(y, z))=U_{T_e}(U_{T_e}(x, y),z)$.
\end{itemize}
\item[3.2.] $y\in (e, 1]$.
\begin{itemize}
\item[3.2.1.] If $z\in [0, e)$, then
\begin{align*}
U_{T_e}(x, U_{T_e}(y, z))&=U_{T_e}(U_{T_e}(z, y), x)~~ (\text{commutativity of } U_{T_e})\\
&= U_{T_e}(z, U_{T_e}(y, x))~~ (\text{by 1.2.3})\\
&= U_{T_e}(U_{T_e}(x, y), z)~~ (\text{commutativity of } U_{T_e}).
\end{align*}

\item[3.2.2.] If $z\in (e, 1]$, then
\begin{align*}
U_{T_e}(x, U_{T_e}(y, z))&= U_{T_e}(U_{T_e}(z, y), x)~~ (\text{commutativity of } U_{T_e})\\
&= U_{T_e}(z, U_{T_e}(y, x))~~ (\text{by 2.2.3})\\
&= U_{T_e}(U_{T_e}(x, y), z)~~ (\text{commutativity of } U_{T_e}).
\end{align*}

\item[3.2.3.] If $z\in I_{e}$, then $U_{T_e}(x, U_{T_e}(y, z))=U_{T_e}(x, z)=x\wedge z$ and
$U_{T_e}(U_{T_e}(x, y), z)=U_{T_e}(x, z)=x\wedge z$, implying that $U_{T_e}(x, U_{T_e}(y, z))=U_{T_e}(U_{T_e}(x, y),z)$.
\end{itemize}

\item[3.3.] $y\in I_{e}$.
\begin{itemize}
\item[3.3.1.] If $z\in [0, e)$, then
\begin{align*}
U_{T_e}(x, U_{T_e}(y, z))&=U_{T_e}(U_{T_e}(z, y), x)~~ (\text{commutativity of } U_{T_e})\\
&= U_{T_e}(z, U_{T_e}(y, x))~~ (\text{by 1.3.3})\\
&= U_{T_e}(U_{T_e}(x, y), z)~~ (\text{commutativity of } U_{T_e}).
\end{align*}

\item[3.3.2.] If $z\in (e, 1]$, then
\begin{align*}
U_{T_e}(x, U_{T_e}(y, z))&= U_{T_e}(U_{T_e}(z, y), x)~~ (\text{commutativity of } U_{T_e})\\
&= U_{T_e}(z, U_{T_e}(y, x))~~ (\text{by 2.3.3})\\
&= U_{T_e}(U_{T_e}(x, y), z)~~ (\text{commutativity of } U_{T_e}).
\end{align*}

\item[3.3.3.] If $z\in I_{e}$, then
\begin{align*}
U_{T_e}(x, U_{T_e}(y, z))&=U_{T_e}(x, y\vee z)\\
&=\begin{cases}
x\wedge y\wedge z, & y\wedge z\in I_{e}, \\
0, & y\wedge z=0,
\end{cases}\\
&=x\wedge y\wedge z,
\end{align*}
and
\begin{align*}
U_{T_e}(U_{T_e}(x, y), z)&=U_{T_e}(x\wedge y, z)\\
&=\begin{cases}
x\wedge y\wedge z, & x\wedge y\in I_{e}, \\
0, & x\wedge y=0,
\end{cases}\\
&=x\wedge y\wedge z,
\end{align*}
implying that $U_{T_e}(x, U_{T_e}(y, z))=U_{T_e}(U_{T_e}(x, y),z)$.
\end{itemize}
\end{itemize}
\end{itemize}
\end{itemize}
\epf

\begin{theorem}\label{U(Se)-operation-Thm}
Let $(X, \leq, 0, 1)$ be a bounded lattice, $e\in L\setminus \{0, 1\}$, and $S_{e}$ be a $t$-conorm on
$[e, 1]$. Then, the following statements are equivalent:
\begin{enumerate}[{\rm (i)}]
\item\label{2.i} $U_{S_{e}}$ defined by \eqref{U(Se)-operation} is a uninorm on $L$ with the neutral
element $e$;

\item\label{2.ii} $\mathscr{P}:=\{x\in (0, e): \exists y\in I_{e}, ~ x\leq y\}=\emptyset$;

\item\label{2.iii} $I_{e}\| (0, e]$.
\end{enumerate}
%%$I_{e} \vee I_{e} \subset I_{e}$ or $I_{e}\vee I_{e}=\{1\}$.
\end{theorem}

\pf
Applying Theorem~\ref{AM-Thm-7}, it is clear that (\ref{2.iii})$\Longleftrightarrow$
(\ref{2.ii})$\Longrightarrow$(\ref{2.i}).

\medskip

(\ref{2.i})$\Longrightarrow$(\ref{2.iii}). Suppose that $I_{e}$ and $(0, e]$ are comparable,
i.e., there exist $x_{0}\in (0, e]$ and $y_{0}\in I_{e}$ such that $x_{0}\leq y_0$.
Clearly, $x_0\neq e$. Since $U_{S_{e}}$ is increasing, applying \eqref{U(Se)-operation} yields that
$$
0<x_0=U_{S_{e}}(x_0, x_0)\leq U_{S_{e}}(x_0, y_0)=0,
$$
which is a contradiction.
\epf

\begin{remark}
\begin{enumerate}[(1)]
\item Theorems~\ref{CK-Thm-3.9} and \ref{CK-Thm-3.12} are respectively direct corollaries of Theorems~\ref{U(T)-operation-Thm}
and \ref{U(S)-operation-Thm}.

\item Theorem~\ref{AM-Thm-6} is a direct corollary of
Theorem~\ref{U(Te)-operation-Thm}, when $\mathscr{P}=\emptyset$.
\end{enumerate}
\end{remark}

\begin{example}\label{Exa-3}
Consider the complete lattice $L$ with the Hasse diagram shown in Fig.~\ref{Hasse-2} and take $T_{e}: \{0, a, e\}^2
\longrightarrow \{0, a, e\}$ by $T_{e}(x, y)=x\wedge y$ for $(x, y)\in \{0, a, e\}^{2}$. Clearly, $I_{e}=\{b, c\}$.
The binary operation $U_T$ defined by \eqref{U(T)-operation} is given in Table~\ref{tab-a}. Applying
Theorem~\ref{U(T)-operation-Thm} shows that $U_T$ is a uninorm on $L$, since $b\vee c=1$. However,
applying Theorem~\ref{CK-Thm-3.9} can not conclude that $U_{T}$ is a uninorm on $L$.
%The binary operations $U_T$ and $U_t$, which are respectively defined by \eqref{U(T)-operation} and \eqref{U(t)-operation*},
%are given by Table~\ref{tab-a} and Table~\ref{tab-b}, respectively. Applying Theorems~\ref{U(T)-operation-Thm} and
%\ref{U(t)-operation-Thm} yields that both $U_T$ and $U_{T_e}$ are uninorms on $L$, as $b\vee c=1$. However,
%Theorems~\ref{CK-Thm-3.9} and \ref{CK-Thm-3.1} do not work.

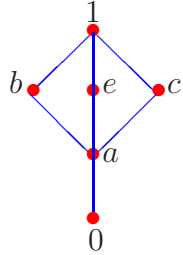
\begin{figure}[h]
\begin{center}
\begin{picture}(150,100)(100,-60)
\put(174.5, 25){$1$}
\put(174.5,19){\color{red}$\bullet$}
\put(176,20){\color{blue}\line(-1,-1){22}}
\put(152,-3.6){\color{red}$\bullet$}
\put(146,-2){$b$}
\put(174.5,-3.6){\color{red}$\bullet$}
\put(181,-2){$e$}
\put(179,20){\color{blue}\line(1,-1){22}}
\put(199,-3.6){\color{red}$\bullet$}
\put(205,-2){$c$}
\put(154,-2.5){\color{blue}\line(1,-1){24}}
\put(201,-2.5){\color{blue}\line(-1,-1){24}}
%%\put(175,-28){\color{blue}\line(0,0){24}}
\put(174.5,-28){\color{red}$\bullet$}
\put(181,-28){$a$}
\put(177.5,-51){\color{blue}\line(0,0){72}}
\put(174.5,-52){\color{red}$\bullet$}
\put(175.5,-61){$0$}
\end{picture}
\end{center}
\caption{The Hasse diagram of the lattice $L$ in Example~\ref{Exa-3}}
\label{Hasse-2}
\end{figure}

\begin{table}[h]
\centering
\begin{tabular}{|l|llllll|}
\hline
 $U_T$& $0$ & $a$ & $b$ & $c$ & $e$ & $1$ \\
 \hline
 ~$0$ & $0$ & $0$ & $b$ & $c$ & $0$ & $1$ \\
 ~$a$ & $0$ & $a$ & $b$ & $c$ & $a$ & $1$ \\
 ~$b$ & $b$ & $b$ & $b$ & $1$ & $b$ & $1$ \\
 ~$c$ & $c$ & $c$ & $1$ & $c$ & $c$ & $1$ \\
 ~$e$ & $0$ & $a$ & $b$ & $c$ & $e$ & $1$ \\
 ~$1$ & $1$ & $1$ & $1$ & $1$ & $1$ & $1$ \\ \hline
\end{tabular}
\label{tab-a}
\caption{The binary operation $U_T$ in Example~\ref{Exa-3}}
\end{table}
\end{example}

\end{document}